\documentclass[final]{acmtrans2m_a}

\usepackage{amsmath}

\usepackage{amsfonts}

\usepackage{amssymb}

\usepackage{stmaryrd}

\usepackage[mathscr]{eucal}

\usepackage{amsbsy}

\usepackage{bm}

\usepackage{paralist}

\usepackage[usenames]{color}

\usepackage[position=below]{subfig}

\usepackage{graphicx}
\usepackage{ifpdf}
\ifpdf
\DeclareGraphicsExtensions{.pdf,.png}
\else
\DeclareGraphicsExtensions{.eps}
\fi

\newcommand{\Sec}[1]{\S\ref{sec:#1}} %section
\newcommand{\Eqn}[1]{{(\ref{eq:#1})}} %equation
\newcommand{\Fig}[1]{{Figure~\ref{fig:#1}}} %figure
\newcommand{\Figs}[2]{{Figures~\ref{fig:#1} and \ref{fig:#2}}} %figure
\newcommand{\Tab}[1]{{Table~\ref{tab:#1}}} %table
\newcommand{\Real}{{\mathbb R}}

\newcommand{\Tra}{^{{\sf T}}} % transpose
\newcommand{\Inv}{^{-1}} % inverse
\newcommand{\V}[1]{{\bm{\mathbf{\MakeLowercase{#1}}}}} % vector
\newcommand{\Vbar}[1]{{\bm \bar{\mathbf{\MakeLowercase{#1}}}}} % vector
\newcommand{\VE}[2]{\V{#1}({#2})} % vector element
\newcommand{\Oprod}{\circ} % outer product

\newcommand{\M}[1]{{\bm{\mathbf{\MakeUppercase{#1}}}}} % matrix
\newcommand{\Mhat}[1]{{\bm{\hat \mathbf{\MakeUppercase{#1}}}}} % matrix
\newcommand{\ME}[3]{\M{#1}(#2,#3)} % matrix element
\newcommand{\MC}[2]{\V{#1}_{#2}} % matrix column
\newcommand{\MhatE}[3]{\Mhat{#1}(#2,#3)} % hatted matrix element
\newcommand{\MbarC}[2]{\Vbar{#1}_{#2}} % matrix column
\newcommand{\T}[1]{\boldsymbol{\mathscr{\MakeUppercase{#1}}}} %tensor
\newcommand{\TE}[4]{\T{#1}(#2,#3,#4)} % tensor element
\newcommand{\TS}[2]{\M{#1}_{{#2}}} % tensor slice
\newcommand{\TSE}[4]{\M{#1}_{{#2}}(#3,#4)} % tensor slice element
\newcommand{\norm}[1]{\left\lVert \, #1 \, \right\rVert}
\newcommand{\nnz}[1]{{\rm nnz}(#1)}

\newcommand{\TheTitle}{Temporal Link Prediction using Matrix and Tensor Factorizations}
\newcommand{\TheShortTitle}{Temporal Link Prediction}

\newcommand{\TheShortAuthors}{D. M. Dunlavy, T. G. Kolda, and E. Acar}
\newcommand{\TheKeywords}{link mining, link prediction, evolution,
  tensor factorization, CANDECOMP, PARAFAC}
\newcommand{\TheAbstract}{%
  The data in many disciplines such as social networks, web
  analysis, etc.\@ is link-based, and the link structure can be
  exploited for many different data mining tasks. In this paper, we
  consider the problem of temporal link prediction: Given link data
  for times 1 through $T$, can we predict the links at time $T+1$?
  If our data has underlying periodic structure, can we predict out
  even further in time, i.e., links at time $T+2$, $T+3$, etc.?
  In this paper, we consider bipartite graphs that evolve over time
  and consider matrix- and tensor-based methods for predicting future links.
  We present a
  weight-based method for collapsing multi-year data into a single
  matrix. We show how the well-known Katz method for link prediction
  can be extended to
  bipartite graphs and, moreover, approximated in a scalable way using a
  truncated singular value decomposition. Using a CANDECOMP/PARAFAC
  tensor decomposition of the data, we illustrate the usefulness of
  exploiting the natural three-dimensional structure of temporal link data.
  Through several numerical experiments, we demonstrate that both matrix-
  and tensor-based techniques are effective for temporal link
  prediction despite the inherent difficulty of the
  problem. Additionally, we show that tensor-based techniques are
  particularly effective for temporal data with varying periodic patterns.
}

  \title{\TheTitle}
    \author{%
      DANIEL M. DUNLAVY and TAMARA G. KOLDA\\
      Sandia National Laboratories \and
      EVRIM ACAR\\
      National Research Institute of Electronics and Cryptology (TUBITAK-UEKAE)
    }
    \begin{abstract}
      \TheAbstract
    \end{abstract}




\begin{document}

\begin{bottomstuff}
  Author Addresses: %
  D.~Dunlavy, Sandia National Laboratories, Albuquerque, NM
  87123-1318, \texttt{dmdunla@sandia.gov}. %
  T.~Kolda, Sandia National Laboratories, Livermore, CA 94551-9159,
  \texttt{tgkolda@sandia.gov}. %
   E.~Acar, TUBITAK-UEKAE, Gebze, Turkey \texttt{evrim.acar@bte.tubitak.gov.tr}.%
\end{bottomstuff}
\maketitle

\section{Introduction}
\label{sec:intro}

The data in different analysis applications such as social networks, communication networks, web
analysis, and collaborative filtering consists of relationships,
which can be considered as links, between objects. For instance, two
people may be linked to each other if they exchange emails or phone
calls.  These relationships can be modeled as a graph, where nodes
correspond to the data objects (e.g., people) and edges correspond to
the links (e.g., a phone call was made between two people).
The link structure of the resulting
graph can be exploited to detect underlying groups of objects,
predict missing links, rank objects, and handle many other tasks
\cite{GeDi05}.

Dynamic interactions over time introduce another dimension to the
challenge of mining and predicting link structure. Here we consider the task of link
prediction in time.  Given link data for $T$ time steps, can we
predict the relationships at time $T+1$? This problem has been
considered in a variety of contexts
\cite{Mohammad2006,LiKl07,Purnamrita2007}.
Collaborative filtering is also a related task, where the objective is
to predict interest of users to objects (movies, books, music) based
on the interests of similar users \cite{Liu2007,Koren2009}.
The \emph{temporal link prediction problem} is different from \emph{missing link
prediction}, which has no temporal aspect and where the goal is to
predict missing connections in order to describe a more complete
picture of the overall link structure in the data \cite{ClMoNe08}.

We extend the problem of
temporal link prediction stated above to the problem of
\emph{periodic} temporal link prediction. For such problems, given
link data for $T$ time steps, can we predict the relationships at
times $T+1, T+2, \dots, T+L$, where $L$ is the length of the periodic
pattern? Such problems often arise in communication networks, such as
e-mail and network traffic data, where weekly or monthly interaction
patterns abound. If we can discover a temporal pattern in the data,
temporal forecasting methods such as Holt-Winters \cite{ChYa88} can be
used to make predictions further out in time.

Time-evolving link data can be organized as a third-order tensor, or multi-dimensional array. In
the simplest case, we can define a tensor $\T{Z}$ of size $M
\times N \times T$ such that
\begin{displaymath}
  \TE{Z}{i}{j}{t} =
  \begin{cases}
    1 & \text{if object $i$ links to object $j$ at time $t$}, \\
    0 & \text{otherwise}.
  \end{cases}
\end{displaymath}
It is also possible to use weights to indicate the strength of the links.
The goal is to predict the links at time $T+1$ or for a period time $T+1, \dots, T+L$ by
analyzing the link structure of $\T{Z}$.

\begin{figure}
  \centering
  \includegraphics[trim=0 5 0 5]{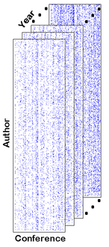}
  \caption{DBLP Data for 1991--2000.}
  \label{fig:data}
\end{figure}

\Fig{data} presents an illustration of such temporal link data for the
1991--2000 DBLP bibliometric data set,
which contains publication data for a large number of
professional conferences in areas related to computer science
(described in more detail in \Sec{numerical}). The
plot depicts the patterns of links between authors and conferences
over time, with blue dots denoting the links (i.e., values of 1 as
defined above).

We consider both matrix- and tensor-based
methods for link prediction. For the matrix-based methods, we collapse
the data into a single matrix by summing (with and without weights) the
matrices corresponding to the time slices. As a baseline, we consider a
low-rank approximation as produced by a truncated singular value decomposition
(TSVD). Next, we consider the Katz method \cite{Ka53}
(extended to bipartite graphs),
which has proven to be highly accurate in previous work on link prediction
\cite{Huang2005,LiKl07,Huang2009}; however, it is not always clear how
to make the method
scalable. Therefore, we present a novel scalable technique for computing
approximate Katz scores based on a truncated spectral decomposition (TKatz).
For the tensor-based methods, we consider the CANDECOMP/PARAFAC (CP)
tensor decomposition \cite{CaCh70,Ha70}, which does not collapse the
data but instead retains its natural three-dimensional structure.
Tensor factorizations are higher-order extensions of matrix
factorizations that capture the underlying patterns
in multi-way data sets and have proved to be successful in diverse
disciplines including chemometrics, neuroscience and social network analysis
\cite{AcYe09,KoBa09}. Moreover, CP yields a highly
interpretable factorization that includes a time dimension. In terms of prediction,
the matrix-based methods are limited to temporal prediction for a single time
step, whereas CP can be used in solving both single step and periodic temporal
link prediction problems.

There are many possible applications for link prediction, such as
predicting the web pages a web surfer may visit on a given day based on
past browsing history, the places that a traveler may fly to in a
given month, or the patterns of computer network traffic.
We consider two applications for link prediction.
First, we consider computer science
conference publication data with a goal of predicting which authors
will publish at which  conferences in year $T+1$ given the publication
data for the previous $T$ years. In this case, we assume we have $M$
authors and $N$ conferences.
All of the methods produce scores for each $(i,j)$ author-conference pair
for a total of $MN$ prediction scores for year $T+1$. For large values of $M$
or $N$, computing all possible
scores is impractical due to the large memory requirements of storing
all $MN$ scores. However, we note that it is possible to easily compute
subsets of the scores. For example, these methods can answer specific
questions such as ``Who is most likely to publish at the KDD conference next year?''
or ``Where is Christos Faloutsos most likely to publish next year?'' using
only $O(M+N)$ memory. This is how we envision link prediction methods
being used in practice.
Second, we consider the problem of how to predict links when
periodic patterns exist in the data. For example, we consider a simulated
example where data is taken daily over ten weeks. We should be able to
recognize, for example, weekday versus weekend patterns and use those
in making predictions. If we consider a scenario of users accessing
various online services, we should be able to differentiate between
services that are heavily accessed on weekdays (e.g., corporate email)
versus those that are used mostly on weekends (e.g., entertainment services).

\subsection{Our Contributions}
The main contributions of this paper can be
summarized as follows:
\begin{compactitem}[$\bullet$]
\item Weighted methods for collapsing temporal data into a matrix are
  shown to outperform straight summation (inspired by the results in \cite{ShNe08})
  in the case of single step temporal link prediction.
\item The Katz method is extended to the case of bipartite graphs and
  its relationship to the matrix SVD is derived.
  Additionally, using the truncated SVD, we devise a scalable method for
  calculating a ``truncated'' Katz score.
\item The CP tensor decomposition is applied to temporal data. We
  provide both heuristic- and forecasting-based prediction methods
  that use the temporal information extracted by CP.
\item Matrix- and tensor-based methods are compared on
  DBLP bibliometric data in terms of link prediction performance
  and relative expense.
\item Tensor-based methods are applied to periodic temporal data with
  multiple period patterns. Using a forecasting-based prediction
  method, it is shown how the method can be used to predict forward in
  time.
\end{compactitem}

\subsection{Notation}
Scalars are denoted by lowercase letters, e.g., $a$. Vectors are
denoted by boldface lowercase letters, e.g., $\V{a}$.
Matrices are denoted by boldface capital letters, e.g., $\M{A}$. The
$r$th column of a matrix $\M{A}$ is denoted by $\MC{A}{r}$.
Higher-order tensors are denoted by boldface Euler script letters,
e.g., $\T{Z}$. The $t$th frontal slice of a tensor $\T{Z}$ is denoted
$\TS{Z}{t}$. The $i$th entry of a vector $\V{a}$ is denoted by $\VE{a}{i}$, element
$(i,j)$ of a matrix $\M{A}$ is denoted by $\ME{A}{i}{j}$, and element $(i,j,k)$
of a third-order tensor $\T{Z}$ is denoted by $\TE{Z}{i}{j}{k}$.

\subsection{Organization}

The organization of this paper is as follows.
Matrix techniques are presented in \Sec{matrix},
including a weighted method for collapsing the
tensor into matrix in \Sec{collapse}, the TSVD method in \Sec{tsvd}, and the Katz and
TKatz methods in \Sec{katz}. The CP tensor technique is presented in
\Sec{tensor}. Numerical results on the DBLP data set are discussed in
\Sec{numerical} and on simulated periodic data in \Sec{numerical2}.
We discuss related work
in \Sec{related}.
Conclusions and future work are discussed in
\Sec{conclusions}.

%%% Local Variables:
%%% mode: latex
%%% TeX-master: "paper"
%%% TeX-command-default: "PDFLaTeX"
%%% End:   % Introduction
\section{Matrix Techniques}
\label{sec:matrix}

We consider different matrix techniques by collapsing the matrices over time into a single matrix.
In \Sec{collapse}, we present two techniques (unweighted and weighted)
for combining the multi-year data into a single matrix. In \Sec{tsvd},
we present the technique of using a truncated SVD to generate link
scores. In \Sec{katz}, we extend the Katz method to bipartite graphs
and show how it can be computed efficiently using a low-rank
approximation.

\subsection{Collapsing the data}
\label{sec:collapse}

Suppose that our data set consists of matrices $\TS{Z}{1}$ through $\TS{Z}{T}$
of size $M \times N$ and the goal is to predict $\TS{Z}{T+1}$.
The most straightforward way to collapse that data into a
single $M \times N$ matrix $\M{X}$ is to sum all the entries across
time, i.e.,
\begin{equation}
  \label{eq:CT}
  \ME{X}{i}{j} = \sum_{t=1}^T \TSE{Z}{t}{i}{j}.
\end{equation}
We call this the \emph{collapsed tensor (CT)} because it collapses
(via a sum) the entries of the tensor
$\T{Z}$ along the time mode.  This is similar to the approach
in \cite{LiKl07}.

\begin{figure}[htbp]
  \centering
  \includegraphics[width=4in,trim=8 5 0 0]{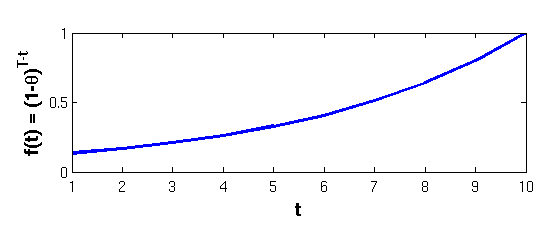}
  \caption{Plot of the decay function $f(t) = (1-\theta)^{T-t}$ for
  $\theta=0.2$ and $T=10$.}
  \label{fig:decay}
\end{figure}

We propose an alternative approach to collapsing the tensor data,
motivated by \cite{ShNe08}, where the
link structure is damped backward in time according to the following
formula:
\begin{equation}
  \label{eq:CWT}
  \ME{X}{i}{j} = \sum_{t=1}^T (1-\theta)^{T-t} \; \TSE{Z}{t}{i}{j}.
\end{equation}
The parameter $\theta \in (0,1)$ can be chosen by the user or
according to experiments on various training data sets.
We call this the \emph{collapsed weighted
  tensor (CWT)} because the slices in the time mode are weighted in
the sum. This gives greater weight to more recent links.
See \Fig{decay} for a plot of $f(t) = (1-\theta)^{T-t}$ for $\theta
= 0.2$ and $T=10$.

The numerical results in \Sec{numerical} demonstrate improved performance using CWT versus CT.

\subsection{Truncated SVD}
\label{sec:tsvd}

\begin{figure}[htbp]
  \centering
  \includegraphics[width=3in]{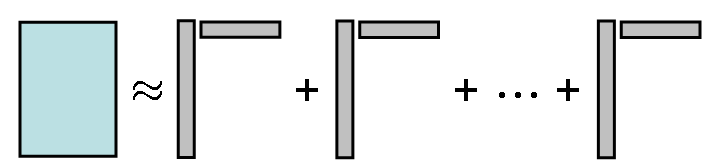}
  \caption{Illustration of the matrix TSVD.%
  }
  \label{fig:svd}
\end{figure}

One of the methods compared in this paper is a
low-rank approximation of the matrix $\M{X}$ produced by \Eqn{CT} or
\Eqn{CWT}. Specifically, suppose that the compact SVD of $\M{X}$ is
given by
\begin{equation}
  \label{eq:svd}
  \M{X} = \M{U} \bm{\Sigma} \M{V}\Tra,
\end{equation}
where $R$ is the rank of $\M{X}$, $\M{U}$ and $\M{V}$ are orthogonal
matrices of sizes $M \times R$ and $N \times R$, respectively, and
$\bm{\Sigma}$ is a diagonal matrix of singular values $\sigma_1 >
\sigma_2 > \cdots > \sigma_R > 0$. It is well known that the best
rank-$K$ approximation of $\M{X}$ is then given by the truncated SVD
\begin{equation}
  \label{eq:tsvd}
  \M{X} \approx \M{U}_K \bm{\Sigma}_K \M{V}_K,
\end{equation}
where $\M{U}_K$ and $\M{V}_K$ comprise the first $K$ columns of
$\M{U}$ and $\M{V}$ and $\bm{\Sigma}_K$ is the $K \times K$ principal
submatrix of $\bm{\Sigma}$.
We can write \Eqn{tsvd}
as a sum of $K$ rank-1 matrices:
\begin{displaymath}
  \M{X} \approx \sum_{k=1}^K \sigma_k \MC{U}{k} \MC{V}{k}\Tra,
\end{displaymath}
where $\MC{U}{k}$ and $\MC{V}{k}$ are the $k$th columns of $\M{U}$ and
$\M{V}$ respectively.
The TSVD is visualized in \Fig{svd}.

A matrix of scores for predicting future links can then be calculated
as
\begin{equation}
  \label{eq:score-tsvd}
  \M{S} = \M{U}_K \bm{\Sigma}_K \M{V}_K.
\end{equation}
We call these the \emph{Truncated SVD (TSVD)} scores.
Low-rank approximations based on the matrix SVD have  proven to be an
effective technique in many data applications; latent semantic
indexing \cite{DuFuLaDe88} is one such example.
This technique is called ``low-rank approximation: matrix entry''
in \cite{LiKl07}.

\subsection{Katz}
\label{sec:katz}

The Katz measure \cite{Ka53} is arguably one of the best link
predictors available because it has been shown to outperform many other
methods \cite{LiKl07}.  Suppose that we have an
undirected graph $G(V,E)$ on $P=|V|$ nodes.  Then the Katz score of a
potential link between nodes $i$ and $j$ is given by
\begin{equation}
\label{eq:katz}
  \MhatE{S}{i}{j} = \sum_{\ell =1}^{+\infty} \beta^{\ell}
  |\text{path}^{(\ell)}_{i,j}|,
\end{equation}
where $|\text{path}^{(\ell)}_{i,j}|$ is the number of paths of length
$\ell$ between nodes $i$ and $j$, and $\beta \in (0,1)$ is a user-defined parameter 
controlling the extent to which longer paths are penalized.

The Katz scores for all pairs of nodes can be expressed in matrix terms as follows. Let
$\Mhat{X}$ be the $P \times P$ symmetric adjacency matrix of the graph.
Then the scores are given by
\begin{equation}
  \label{eq:symkatz}
  \Mhat{S} = \sum_{\ell=1}^{+\infty} \beta^{\ell} \Mhat{X}^\ell
  = (\M{I} - \beta \Mhat{X})\Inv - \M{I}.
\end{equation}
Here $\M{I}$ is the $P \times P$ identity matrix.
If the graph under consideration has weighted edges, $\Mhat{X}$ is
replaced by a weighted adjacency matrix.

We address two problems with the formulation of the Katz measure.
First, the method is not scalable because it requires the inversion of
a $P \times P$ matrix at a cost of $O(P^3)$ operations. We shall see that
we can replace $\Mhat{X}$ by a low-rank approximation in order to compute
the Katz scores more efficiently.
Second, the method is only applicable to square symmetric matrices
representing undirected graphs. We show that it can also be applied to
our situation: a rectangular matrix representing a bipartite graph.

\subsubsection{Truncated Katz}
Assume $\Mhat{X}$ has rank $R \leq P$.
Let the eigendecomposition of $\Mhat{X}$ be given by
\begin{equation}
  \label{eq:evd}
  \Mhat{X} = \Mhat{W} \bm{\hat \Lambda} \Mhat{W}\Tra,
\end{equation}
where $\Mhat{W}$ is a $P \times P$ orthogonal matrix%
\footnote{Recall that if $\M{W}$ is an orthogonal matrix $\M{W}$, then
  $\M{W}\M{W}\Tra = \M{W}\Tra\M{W} = \M{I}$, $\M{W}\Inv = \M{W}\Tra$,
  and $(\M{W}\Tra)\Inv = \M{W}$.} %
and $\bm{\hat\Lambda}$ is a diagonal matrix with $|\hat\lambda_1| \geq
|\hat\lambda_2| \geq \cdots \geq |\hat\lambda_R| > \hat\lambda_{R+1} = \cdots =
\hat\lambda_P = 0$. Then the Katz scores in
\Eqn{symkatz} become
\begin{align*}
  \Mhat{S}
  &  = (\M{I} - \beta \Mhat{W}\bm{\hat\Lambda}\Mhat{W}\Tra)\Inv - \M{I} \\
 &  = \Mhat{W}  \left[ (\M{I} - \beta \bm{\hat\Lambda})\Inv - \M{I} \right]
  \Mhat{W}\Tra
   = \Mhat{W} \bm{\hat \Gamma} \Mhat{W}\Tra,
\end{align*}
where $\bm{\hat \Gamma}$ is a $P \times P$ diagonal matrix with diagonal entries
\begin{equation}
  \label{eq:Katz_scaling}
  \hat \gamma_p = \frac{1}{1 - \beta \hat\lambda_p} - 1
  \quad \text{for} \quad
  p = 1,\dots, P.
\end{equation}
Observe that $\hat \gamma_p=0$ for $p > R$. Therefore,
without loss of generality, we can assume that $\Mhat{W}$ and
$\bm{\hat \Gamma}$ are given in compact form, i.e., $\bm{\hat \Gamma}$ is just
an $R \times R$ diagonal matrix and $\Mhat{W}$ is a $P \times R$
orthogonal matrix.

This shows a close relationship between the Katz measure and the
  eigendecomposition and gives some hint as to how to incorporate a
  low-rank approximation.
The best rank-$L$ approximation of $\Mhat{X}$ is given by replacing $\bm{
  \hat \Lambda}$ in \Eqn{evd} with a matrix $\bm{\hat\Lambda}_L$ where all but
the $L$ largest magnitude diagonal entries are set to zero. The
mathematics carries through as above, and the end result is that the
  Katz scores based on the rank-$L$ approximation are
\begin{displaymath}
  \Mhat{S} = \Mhat{W}_L \bm{\hat\Gamma}_L \Mhat{W}_L\Tra
\end{displaymath}
where $\bm{\hat\Gamma}_L$ is the $L \times L$ principal submatrix of
$\bm{\hat\Gamma}$, and $\Mhat{W}_L$ is the $P \times L$ matrix containing the first $L$ columns of $\Mhat{W}$.

Since it is possible to construct a rank-L approximation of the
adjacency matrix in $O(L|E|)$ operations (using an Arnoldi or Lanczos
technique \cite{Sa92}), this technique can be applied to large-scale
problems at a relatively low cost. We note that in \cite{LiKl07}, Katz
is applied to a low-rank approximation of the adjacency matrix which
is equivalent to what we discuss here, but
its computation is not discussed --- specifically, the fact that it
can be computed efficiently via the formula above is not mentioned. Thus, we assume
that calculation was done directly on the dense low-rank approximation
matrix given by
\begin{displaymath}
  \Mhat{X}_L = \Mhat{W}_L \bm{\hat\Lambda}_L \Mhat{W}_L\Tra.
\end{displaymath}

We contrast this with the approach of Wang et al. \citeyear{Wang2007} who discuss
an approximate Katz measure given by truncating the sum in \Eqn{katz}
to the first $L$ terms (they recommend $L=4$), i.e., $\Mhat{S} =
\sum_{\ell=1}^4 \beta^{\ell} \Mhat{X}^{\ell}$; the main drawback of
this approach is the power matrices may be dense, depending on the
connectivity of the graph.

\subsubsection{Bipartite Katz \& Truncated Bipartite Katz}

Our problem is different than what we have discussed so far because we
are considering a bipartite graph,
represented by a weighted adjacency matrix from \Eqn{CT} or
\Eqn{CWT}. This can be considered as a graph on $P=M+N$ nodes
where the weighted adjacency matrix is given by
\begin{displaymath}
  \Mhat{X} =
  \begin{bmatrix}
    \M{0} & \M{X} \\
    \M{X}\Tra & \M{0}
  \end{bmatrix}.
\end{displaymath}
If $\M{X}$ is rank $R$ and its SVD is given as in $\Eqn{svd}$, then
the eigenvectors and eigenvalues of $\Mhat{X}$ are given by
\def\foo{\frac{1}{\sqrt{2}}}
\begin{displaymath}
  \Mhat{W} =
  \begin{bmatrix}
    \foo \M{U} & - \foo \M{U} \\
    \foo \M{V} & \foo \M{V}
  \end{bmatrix}
  \quad \text{and} \quad
  \bm{\hat\Lambda} =
  \begin{bmatrix}
    \bm{\Sigma} & \M{0} \\
    \M{0} & -\bm{\Sigma}
  \end{bmatrix}.
\end{displaymath}
Note that the eigenvalues in $\bm{\hat \Lambda}$ are not sorted by magnitude
and the rank of $\Mhat{X}$ is $2R$.
Scaling the eigenvalues in $\Mhat{\Lambda}$ as in \Eqn{Katz_scaling}, we get the $2R \times 2R$ diagonal matrix $\Mhat{\Gamma}$ with entries \begin{equation*}
  \hat{\gamma_p} = \left\{ 
  \begin{array}{ll}
  \displaystyle\frac{1}{1 - \beta \sigma_p} - 1 &\quad \text{for} \quad p=1,\dots,R \; , \; \mbox{and}\\
  \displaystyle\frac{1}{1 + \beta \sigma_{p-R}} - 1 &\quad \text{for} \quad p=R+1,\dots,2R.
  \end{array} \right.
\end{equation*}
The square matrix of Katz scores is then given by
\begin{displaymath}
  \Mhat{S} = \Mhat{W}\Mhat{\Gamma}\Mhat{W}\Tra = 
  \begin{bmatrix}
    \M{U}\M{\Psi}^+\M{U}\Tra & \M{U}\M{\Psi}^-\M{V}\Tra \\
    \M{V}\M{\Psi}^-\M{U}\Tra & \M{V}\M{\Psi}^+\M{V}\Tra 
  \end{bmatrix} 
\end{displaymath}
where $\M{\Psi}^-$ and $\M{\Psi}^+$ are diagonal matrices with entries
\begin{eqnarray}
  \psi_p^- &= \displaystyle\frac12\left(\left(\displaystyle\frac{1}{1 - \beta \sigma_p} - 1\right) - \left(\displaystyle\frac{1}{1 + \beta \sigma_p} - 1\right)\right) = \displaystyle\frac{\beta\sigma_p}{1 - \beta^2 \sigma_p^2} \; , \; \mbox{and} \label{eq:psi_minus_p}\\
  \psi_p^+ &=\displaystyle\frac12\left( \left(\displaystyle\frac{1}{1 - \beta \sigma_p} - 1\right) + \left(\displaystyle\frac{1}{1 + \beta \sigma_p} - 1\right)\right) =\displaystyle\frac{1}{1 - \beta^2 \sigma_p^2} - 1 , \label{eq:psi_plus_p}
\end{eqnarray}
respectively, for $p=1,\dots,R$. The link scores for the bipartite graph can be extracted and are given
by
\begin{equation}
  \label{eq:katzscore}
  \M{S} = \M{U}\M{\Psi}^-\M{V}\Tra.
\end{equation}
We call these the \emph{Katz} scores.

We can replace $\M{X}$ by its
best rank-$K$ approximation as in \Eqn{tsvd}, and the resulting Katz
scores then become
\begin{equation}
  \label{eq:tkatzscore}
  \M{S} = \M{U}_K\M{\Psi}^-_K\M{V}_K\Tra,
\end{equation}
where $\M{\Psi}^-_K$ is the $K \times K$ principal submatrix of
$\M{\Psi}^-$. We call these the \emph{Truncated Katz (TKatz)} scores. It is
interesting to note that TKatz is very similar to using TSVD
except that the diagonal weights have been changed. Related
methods for scaling have also been proposed in the area of
information retrieval (e.g., \cite{Bast2005,Yan2008}) where exponential
scaling of singular values led to improved performance.

\subsection{Computational Complexity and Memory}

Computing a sparse rank-$K$ TSVD via an Arnoldi or Lanczos method requires
$O(\nnz{\M{X}})$ work per iteration where $\nnz{\M{X}}$ is the number
of nonzeros
in the adjacency matrix $\M{X}$, which is equal to the number of edges
in the bipartite graph. The number of iterations is typically a small
multiple of $K$ but cannot be known in advance. The
storage of the factorization requires only $K(M+N+1)$
space for the singular values and two factor matrices.
Because TKatz is based on the TSVD, it requires the same amount of
computation and storage for a rank-$K$ approximation. The only
difference is that TKatz stores $\M{\Psi}^-_K$ rather than
$\M{\Sigma}_K$. Katz, on the other hand, requires $O(M^2N + MN^2 + N^3)$ operations
to compute \Eqn{symkatz} if
$M>N$. Furthermore, it stores all of the scores explicitly,
using $O(MN)$ storage.

%%% Local Variables:
%%% mode: latex
%%% TeX-master: "paper"
%%% TeX-command-default: "PDFLaTeX"
%%% End: 
\section{Tensor Techniques}
\label{sec:tensor}

The data set consisting of matrices $\TS{Z}{1}$ through $\TS{Z}{T}$ is three-way,
so this lends itself to a multi-dimensional
interpretation. By analyzing this data set using a three-way model, we can
explicitly model the time dimension and have no need to collapse
the data as discussed in \Sec{collapse}.

\subsection{CP Tensor Model}

\begin{figure}[htbp]
  \centering
  \includegraphics[width=3in]{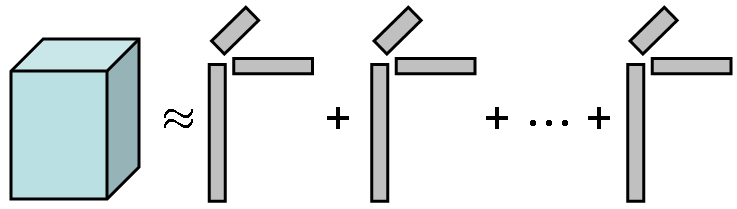}
  \caption{Illustration of the tensor CP model.}
  \label{fig:cp}
\end{figure}

One of the most common and useful tensor models is CP \cite{CaCh70,Ha70};
see also reviews \cite{AcYe09,KoBa09}. Given a three-way tensor $\T{Z}$ of
size $M \times N \times T$, its $K$-component CP decomposition is given by
\begin{equation}
  \label{eq:cp}
  \T{Z} \approx \sum_{k=1}^K \lambda_k \; \MC{A}{k} \Oprod \MC{B}{k} \Oprod \MC{C}{k}.
\end{equation}
Here the symbol $\Oprod$ denotes the outer product%
\footnote{A three way outer product is defined as follows: $\T{X} =
  \V{a} \Oprod \V{b} \Oprod \V{c}$ means
  $\TE{X}{i}{j}{k} = \VE{a}{i} \VE{b}{j} \VE{c}{k}$.}, %
$\lambda_k \in \Real_+$, $\MC{A}{k} \in \Real^M$, $\MC{B}{k} \in
\Real^N$, and $\MC{C}{k} \in \Real^T$ for $k=1,\dots,K$.  Each summand
($\lambda_k \; \MC{A}{k} \Oprod \MC{B}{k} \Oprod \MC{C}{k}$) is called
a \emph{component}, and the individual vectors are called
\emph{factors}.  We assume $\norm{\MC{A}{k}} = \norm{\MC{B}{k}} =
\norm{\MC{C}{k}} = 1$ and therefore $\lambda_k$ contains the scalar
weight of the $k$th component.
An illustration of CP is shown in \Fig{cp}.

The CP tensor decomposition can be considered an analogue of the SVD because it
decomposes a tensor as a sum of rank-one tensors just as the SVD
decomposes a matrix as a sum of rank-one matrices as shown in \Fig{svd}.
Nevertheless, there are also important differences between these decompositions. The
columns of $\M{U}$ and $\M{V}$ are orthogonal in the SVD while there
is no orthogonality constraint in the CP model. Despite the CP model's
lack of orthogonality, Kruskal \cite{Kr89,KoBa09} has shown that CP components are
unique, up to permutation and scaling, under mild conditions. It is because of this property that
we use CP model in our studies. The uniqueness of CP enables the use of factors directly for
forecasting as discussed in \Sec{forecasting}. On the other hand, some other tensor  models such as Tucker \cite{Tucker63,Tu66}
suffers from rotational freedom and factors in the time mode, thus the forecasts, can easily change depending on the rotation
applied to the factors. We leave whether or not such models would be applicable for link prediction as a topic of future research.

\subsection{CP Scoring using a Heuristic}
\label{sec:heuristic}
We make use of the components extracted by the CP model to assign
scores to each pair $(i,j)$ according to their likelihood of
linking in the future.
The outer product of $\MC{A}{k}$ and
$\MC{B}{k}$, i.e., $\MC{A}{k}\MC{B}{k}\Tra$, quantifies the
relationship between object pairs in component $k$.
The temporal profiles are captured in the vectors $\MC{C}{k}$.
Different components may have different trends, e.g.,
they may have
increasing, decreasing, or steady profiles.
In our heuristic approach, we assume that average activity in the last $T_0=3$
years is a good choice for the weight.
We define the similarity score for objects $i$ and $j$ using
a $K$-component CP model in \Eqn{cp} as the $(i,j)$ entry of the following matrix:
\begin{equation}
  \label{eq:score-cp}
    \M{S} = \sum_{k=1}^K \gamma_k\lambda_k\MC{A}{k}\MC{B}{k}\Tra,
    \quad\text{where}\quad
    \gamma_k= \frac{1}{T_0} \sum_{t=T-T_0+1}^T\MC{C}{k}(t).
\end{equation}
This is a simple approach, using temporal information from the
last $T_0=3$ time steps only.
In many cases, the simple heuristic of just averaging the last few
time steps works quite well and is sufficient. An alternative that
provides a more sophisticated use of time is discussed in the next
section.

\subsection{CP Scoring using Temporal Forecasting}
\label{sec:forecasting}
Alternatively, we can use the
temporal profiles computed by CP as a basis for predicting the scores
in future time steps. In this work, we use the Holt-Winters forecasting
method \cite{ChYa88}, which is particularly suitable for time-series
data with periodic patterns. This is an automatic method which only requires the data and the
expected period (e.g., we use $L=7$ for daily data). As will be shown
in \Sec{numerical2}, the Holt-Winters method is fairly accurate in
picking up patterns in time and therefore can be used as a predictive
tool. If we are predicting for $L$ time steps in the future (one period), we get a
tensor of prediction scores of size $M \times N \times L$. This is
computed as
\begin{equation}
  \label{eq:score2-cp}
  \T{S} = \sum_{k=1}^K \lambda_k \MC{A}{k} \Oprod \MC{B}{k} \Oprod
  \bm{\gamma}_k
\end{equation}
where each $\bm{\gamma}_k$ is a vector of length $L$ that is the
prediction for the next $L$ time steps from the Holt-Winters methods
with $\MC{C}{k}$ as input.

\begin{figure}[t!]
  \centering
  %Trim = amount to cut off on LEFT, BOTTOM, RIGHT, TOP
  \includegraphics[trim = 20 15 20 15,clip,height=2.5in]{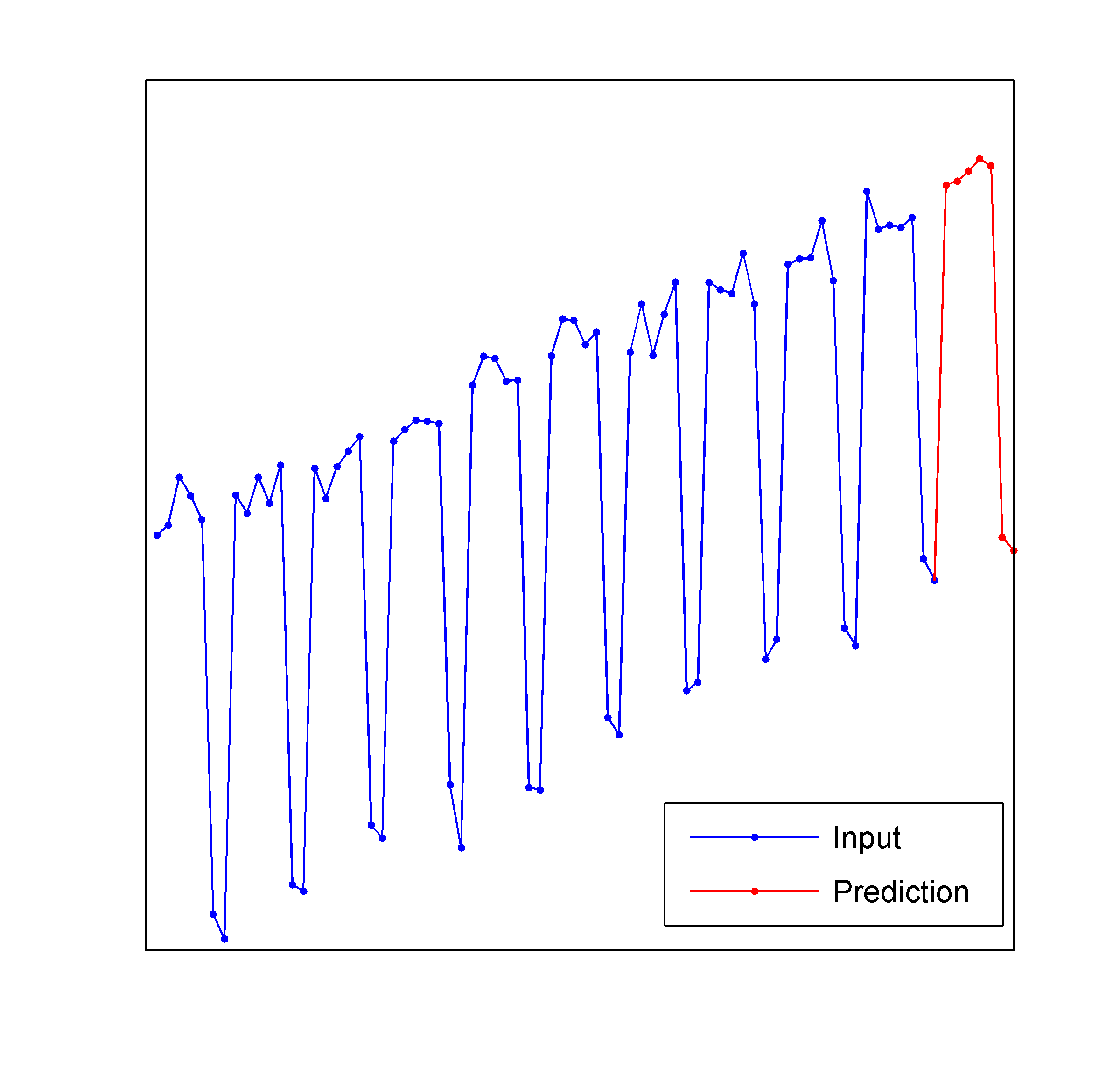}
  \caption{An illustration of the predictions produced by the
    additive Holt-Winters method on data with a period of $L=7$.}
  \label{fig:hw}
\end{figure}

For our studies, we implemented the \emph{additive} Holt-Winters
method (as described in \citeN{ChYa88}), i.e., Holt's linear trend model with additive seasonality,
which corresponds to an exponential smoothing
method with additive trend and additive seasonality. For a review of exponential smoothing methods,
see \cite{Ga06}. An example of forecasting using additive Holt-Winters is
shown in \Fig{hw}. The input is shown in blue, and the prediction of
the next $L=7$ time steps is shown in red. We show examples in
\Sec{numerical2} that use the actual CP data as input.
Forecasting methods beside the Holt-Winters method have also proven useful in analyzing temporal data \cite{MaHi00}.
Work on the applicability of different forecasting methods for link prediction
is left for future work.

\subsection{Computational Complexity and Memory}

The computational complexity of CP is $O(\nnz{\T{Z}})$ per iteration. As with TSVD, we cannot predict the
number of iterations in advance. The storage required for CP is
$K(M+N+T+1)$, for the three factor matrices and the scalar $\lambda_k$
values.

%%% Local Variables:
%%% mode: latex
%%% TeX-master: "paper"
%%% TeX-command-default: "PDFLaTeX"
%%% End: 

\section{Experiments with Link Prediction for One Time Step}
\label{sec:numerical}
We use the DBLP data set%
\footnote{\url{http://www.informatik.uni-trier.de/~ley/db/index.html}} %
to assess the performance of various link
predictors discussed in \Sec{matrix} and \Sec{tensor}. All experiments
were performed using Matlab 7.8 on a Linux Workstation
(RedHat 5.2) with 2 Quad-Core Intel Xeon 3.0GHz processors and 32GB RAM.
We compute the CP model via an Alternating Least Squares (ALS)
approach using the Tensor Toolbox for Matlab \cite{BaKo07}.

\subsection{Data}
At the time the DBLP data was downloaded for this work, it contained publications
from 1936 through the end of 2007. Here we only consider publications of
type \emph{inproceedings} between 1991 and 2007\footnote{The publications
between 1936 and 1990 comprise only 6\% of
publications of type \emph{inproceedings}.}.

The data is organized as a third-order tensor $\T{Z}$ of size $M
\times N \times T$. We let $\TE{C}{i}{j}{t}$ denote the total number
of papers by author $i$ at conference $j$ in year $t$.
In order to decrease the effect of large
numbers of publications, we preprocess the data so that
\begin{equation*}
\TE{Z}{i}{j}{t}=
\begin{cases}
  1 + \log(\TE{C}{i}{j}{t}) & \text{if } \TE{C}{i}{j}{t} > 0,\\
  0 & \text{otherwise}.
\end{cases}
\end{equation*}
Using a sliding window approach, we
divide the data into seven training/test sets such that each training
set contains $T=10$ years and the corresponding test set contains
the following 11th year. \Tab{Data} shows the size and density of the
training and testing sets.
We only keep those authors that have at least 10 publications (i.e., an
average of one per year) in the training data, and each test
set contains only the authors and conferences available in the
corresponding training set.

\begin{table}[t]
  \footnotesize
  \centering
  \setlength{\tabcolsep}{4pt}
  \begin{tabular}{|c c c c c c c|}
    \hline
    \bf Training &  \bf Test & \bf Authors & \bf Confs. & \bf Training Links & \bf Test Links  & \bf  Test New Links  \\
    \bf Years & \bf Year  & &  & \bf (\% Density) & \bf (\% Density)  & \bf (\% Density)  \\
    \hline
    1991-2000 &  2001 & 7108 & 1103 & 112,730 (0.14) & 12,596 (0.16) &
    5,079 (0.06) \\ \hline
    1992-2001 &  2002 & 8368 & 1211 & 134,538 (0.13) & 16,115 (0.16) &
    6,893 (0.07) \\ \hline
    1993-2002 &  2003 & 9929 & 1342 & 162,357 (0.12) & 20,261 (0.15) &
    8,885 (0.07) \\ \hline
    1994-2003 &  2004 & 11836 & 1491 & 196,950 (0.11) & 27,398 (0.16)
    & 12,738 (0.07) \\ \hline
    1995-2004 &  2005 & 14487 & 1654 & 245,380 (0.10) & 35,089 (0.15)
    & 16,980 (0.07) \\ \hline
    1996-2005 &  2006 & 17811 & 1806 & 308,054 (0.10) & 40,237 (0.13)
    & 19,379 (0.06) \\ \hline
    1997-2006 &  2007 & 21328 & 1934 & 377,202 (0.09) & 41,300 (0.10)
    & 20,185 (0.05) \\ \hline
  \end{tabular}
   \caption{Training and Test set pairs formed from the DBLP data set.}
  \label{tab:Data}
\end{table}

\subsection{Interpretation of CP}
\begin{figure}[t!]
  \centering
  \subfloat[Factors from component 3: Top authors are Alberto L. Sangiovanni Vincentelli,
  Robert K. Brayton, Sudhakar M. Reddy, and Irith Pomeranz.
  Top conferences are DAC, ICCAD, and ICCD.]
  {
    %Trim = amount to cut off on LEFT, BOTTOM, RIGHT, TOP
    \includegraphics[width=2.8in,trim=20 8 15 10,clip]{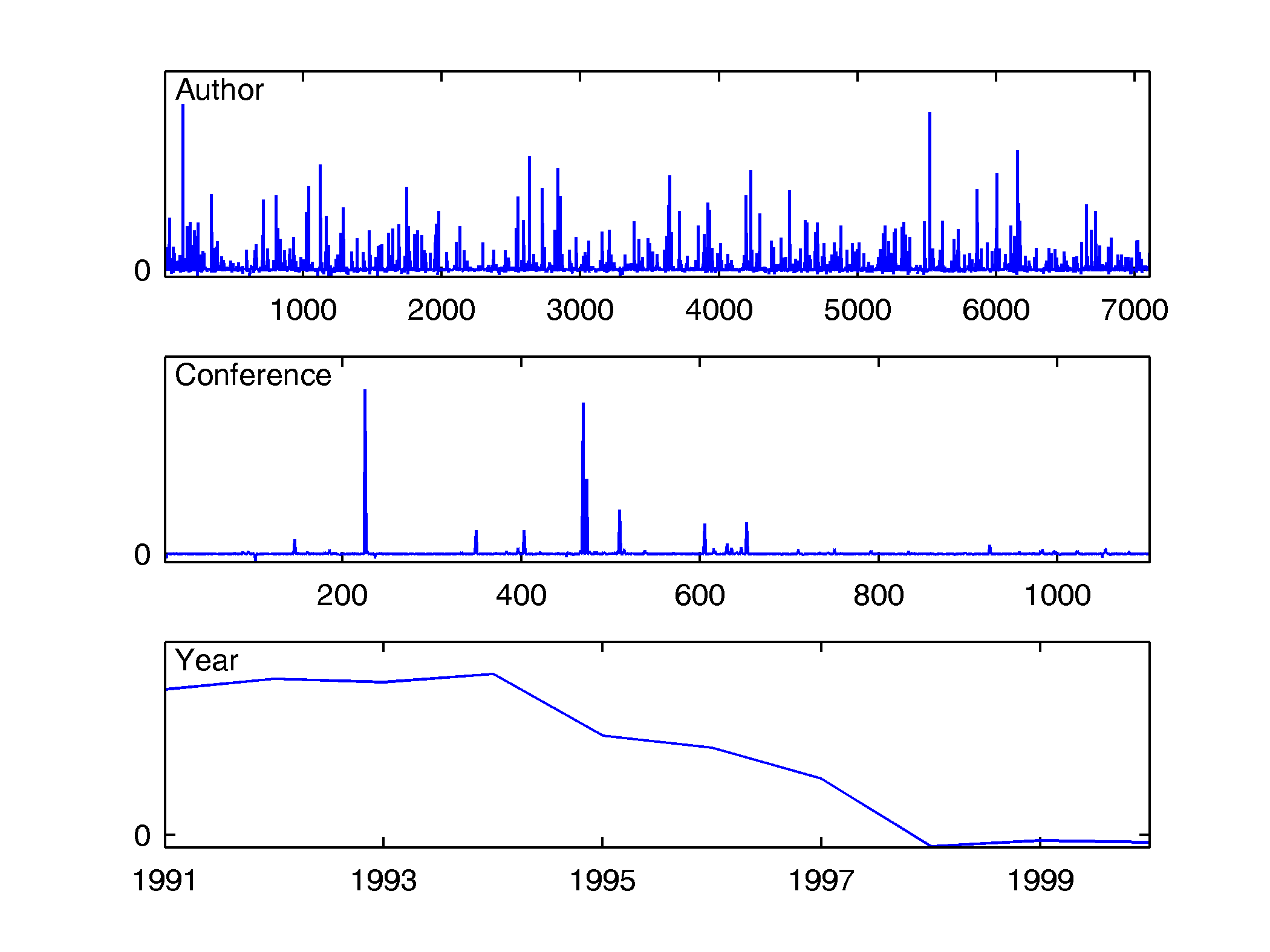}
    \label{fig:c3}
  }\\
  \subfloat[Factors from component 4: Top authors are Miodrag Potkonjak, Massoud Pedram, Jason Cong, and Andrew B. Kahng.
  Top conferences are DAC, ICCAD, and ASPDAC.]
  {
    %Trim = amount to cut off on LEFT, BOTTOM, RIGHT, TOP
    \includegraphics[width=2.8in,trim=20 8 15 10,clip]{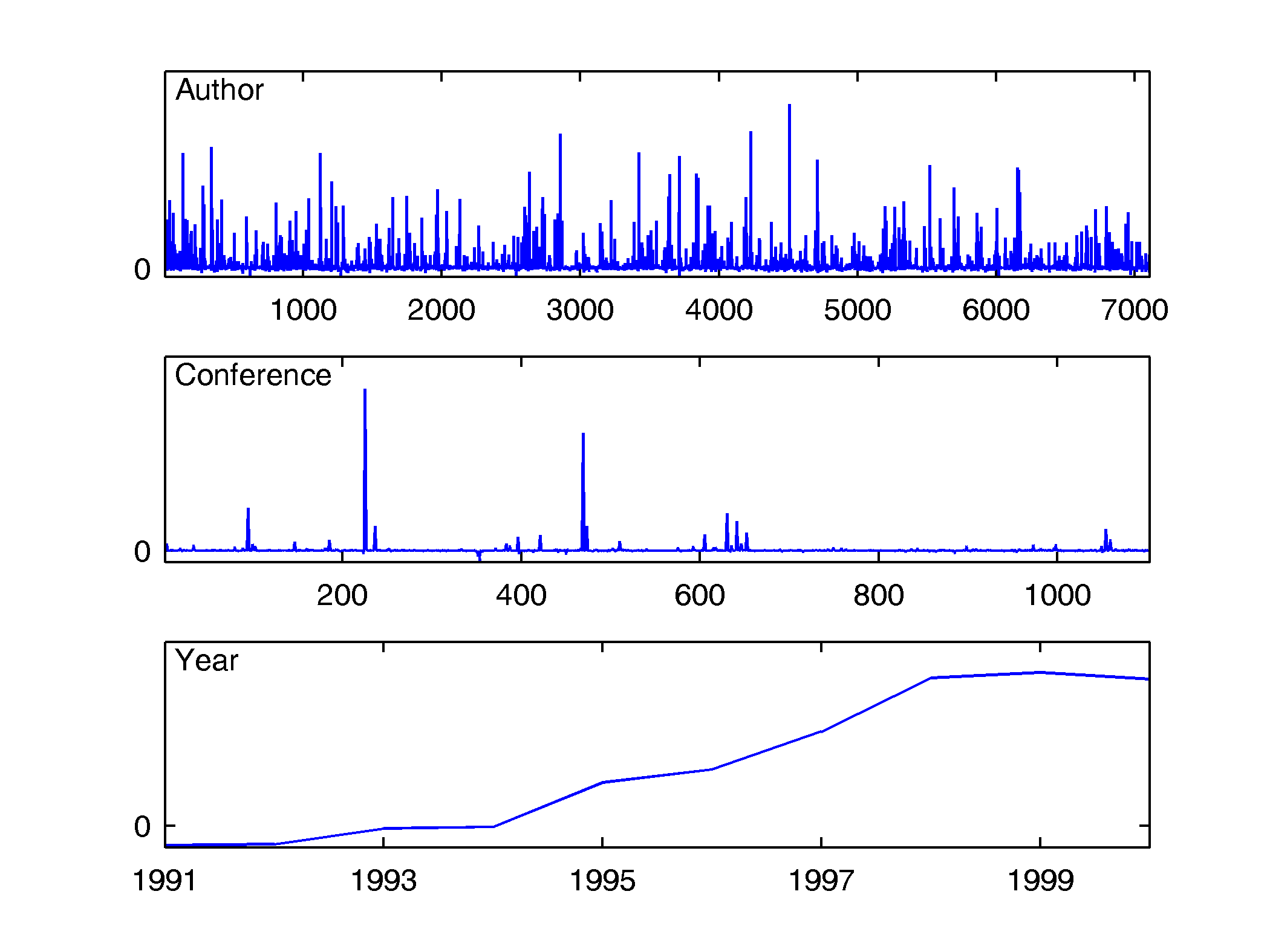}
    \label{fig:c4}
  }
  \caption{Examples from 50-component CP model of publications from 1991 to 2000.}
  \label{fig:components}
\end{figure}

Before addressing the link prediction problem, we first discuss how to use the CP model for exploratory analysis
of the temporal data. The primary advantage of the CP model is its interpretability, as
illustrated in \Fig{components}, which contains three example
components from the 50-component CP model of the tensor representing
publications from 1991 to 2000.
The factor $\MC{A}{k}$ captures a certain group of authors while
$\MC{B}{k}$ extracts the conferences, where the authors captured by
$\MC{A}{k}$ publish. Finally, $\MC{C}{k}$ corresponds to the temporal
signature, depicting the pattern of the publication history
of those authors at those conferences over the associated time period. Therefore, the CP model
can address the link prediction problem well by capturing the
evolution of the links between objects using the factors in the time
mode.

\begin{figure}[htbp]
  \centering
    %Trim = amount to cut off on LEFT, BOTTOM, RIGHT, TOP
    \includegraphics[width=2.8in,trim=20 8 15 10,clip]{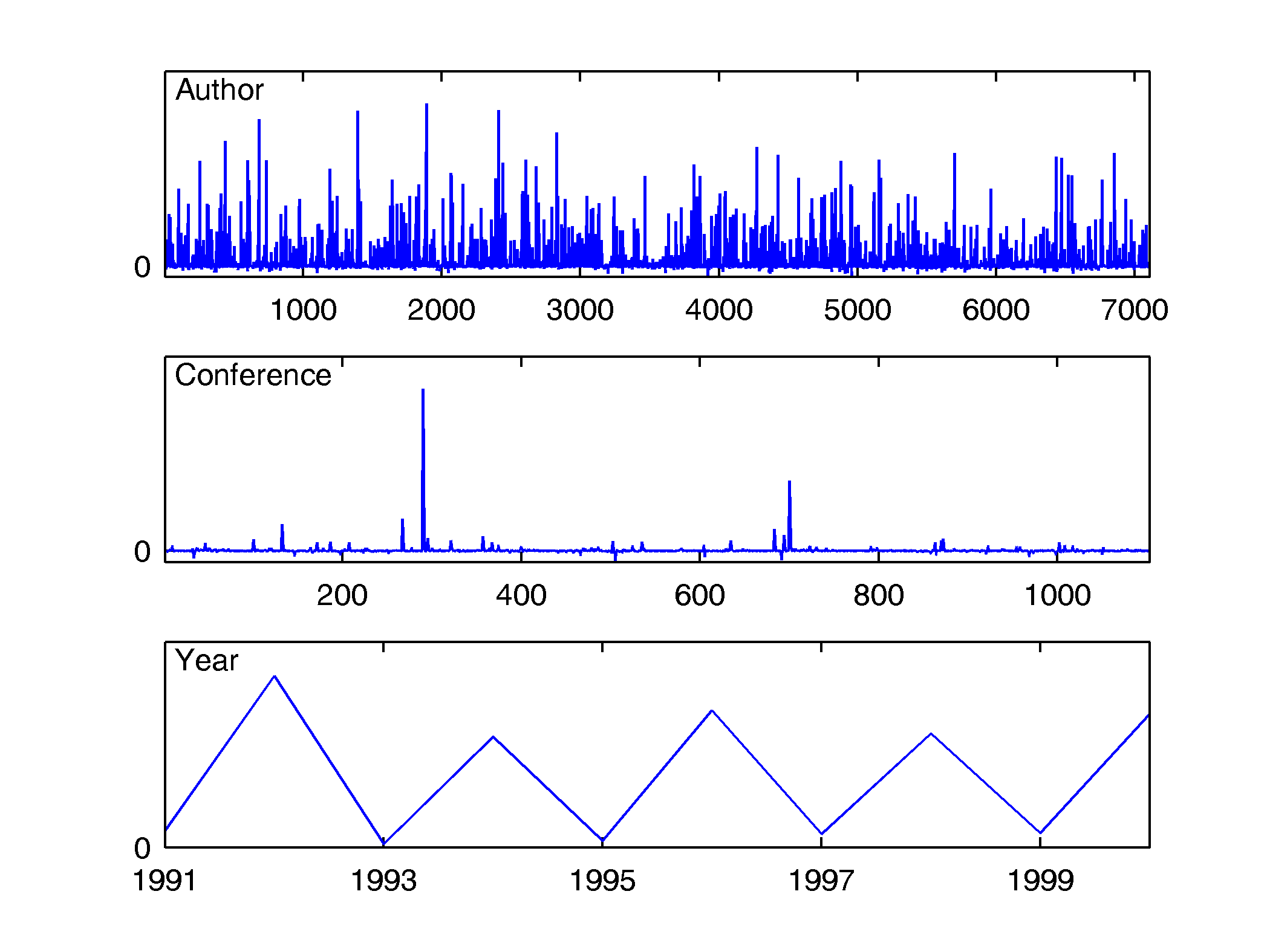}
  \caption{Factors from component 46 of 50-component CP model of
    publications from 1991 to 2000: Top authors are Franz Baader,
    Henri Prade,
  Didier Dubois, and Bernhard Nebel. Top conferences are ECAI and KR.}
    \label{fig:c46}
\end{figure}

\Fig{c3} shows the third component with authors ($\MC{A}{k}$) in the top plot,
conferences ($\MC{B}{k}$) in the middle, and time ($\MC{C}{k}$) on the bottom. The highest scoring
conferences are DAC, ICCAD and ICCD, which are related conferences
on computer design. Many authors publish in these conferences between
1991 and 2000, but the top are Vincentelli, Brayton, and others listed in the caption. This
author/conference combination has a peak in the early 1990s and starts to decline in mid-'90s.
Note that the author and conference scores are mostly positive.
\Fig{c4} shows another example component, which actually has very
similar conferences to those in the component discussed above.
The leading authors, however, are different.
Moreover, the time profile is different with an increasing trend after
the mid-'90s.
\Fig{c46} shows a component that detects
related conferences that take place only in even years. Again we see that
the components are primarily positive. A nice feature of the CP model
is that it does not have any constraints (like orthogonality in the
SVD) that artificially impose a need for negative entries in the
components.

\subsection{Methods and Parameter Selection}

The goal of a link predictor in this study is to predict whether the
$i$th author is going to publish at the $j$th conference during the
test year. Therefore, each nonzero entry in the test set is treated as
1, i.e., a positive link, regardless of the actual number of
publications; otherwise, it is 0 indicating that there is
no link between the corresponding author-conference pair.

The common parameter for all link predictors, except Katz-CT/CWT,
is the number of components, $K$. In our experiments, instead
of using a specific value of $K$, which cannot be determined
systematically, we use an ensemble approach. Let $\TS{S}{K}$ denote the
matrix of scores computed for $K=10,20,...100$. Then the
matrix of ensemble scores, $\M{S}$, used for link prediction is
calculated as
\begin{displaymath}
  \M{S} = \sum_{K \in \{10,20,...100\}} \frac{\TS{S}{K}}{\norm{\TS{S}{K}}_F} \; .
\end{displaymath}
In addition to the number of components, the parameter $\beta$ used in the Katz
scores in \Eqn{katzscore} and \Eqn{tkatzscore} needs to be determined.
We use
$\beta=0.001$, which was chosen such that $\psi^-_p > 0$
for all $p=1,\dots,R$ in \Eqn{psi_minus_p} for the data in our
experiments.
We have observed that if $\psi^-_p < 0$ then
the scores contain entries with large magnitudes but negative
values, which degrades the performance of Katz measure.
Finally, $\theta$ is the parameter used for weighting slices while
forming the CWT in \Eqn{CWT}. We set $\theta = 0.2$
according to preliminary tests on the training data sets.
We use the heuristic scoring method discussed in \Sec{heuristic} for
CP.

\subsection{Link Prediction Results}
\renewcommand{\textfraction}{0.05}
Two experimental set-ups are used to evaluate the
performance of the methods.
\begin{itemize}
\item{\emph{Predicting All Links}: The first approach compares the
    methods in terms of how well they predict positive links in the
    test set}.
\item{\emph{Predicting New Links}: The second approach
    addresses a more challenging problem, i.e., how well the methods
    predict the links that have not been previously seen at any time
    in the training set.}
\end{itemize}

\begin{figure}[t!]
  \centering
  \includegraphics[width=3in]{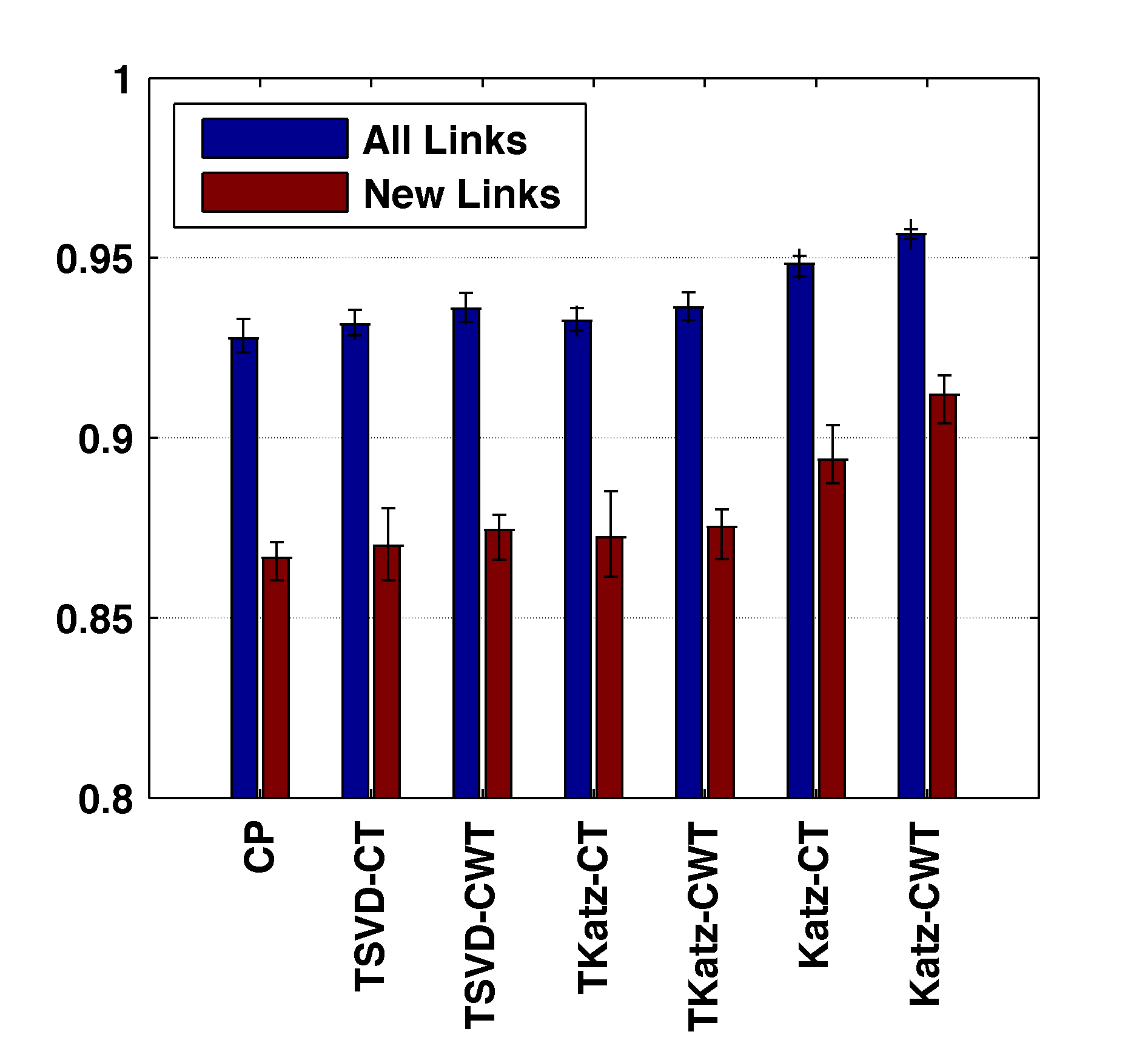}
  \caption{Average link prediction performance of each method across
  all seven training/test set pairs (black bars show absolute range).}
  \label{fig:PerfComp}
\end{figure}

As an evaluation metric for link prediction performance, we use
the area under the receiver operating characteristic curve (AUC)
because it is viewed as a robust measure in the presence of imbalance
\cite{StLuTr06}, which is important since less than 0.2\% of all
possible links exist in our testing data.
\Fig{PerfComp} shows the performance of each link predictor
in terms of AUC when predicting all links (blue bars) and
new links (red bars). As expected, the AUC values are much lower for
the new links. Among all methods, the best performing method in terms
of AUC is Katz-CWT. Further, CWT is consistently better on average than the
corresponding CT methods, which shows that giving more weight to the
data in recent years improves link prediction.

\begin{figure}[t!]
  \centering
  \subfloat[Prediction of \emph{all} links in the test sets.]
  {
    %Trim = amount to cut off on LEFT, BOTTOM, RIGHT, TOP
    \includegraphics[width=.45\textwidth,trim=0 0 0 0,clip]{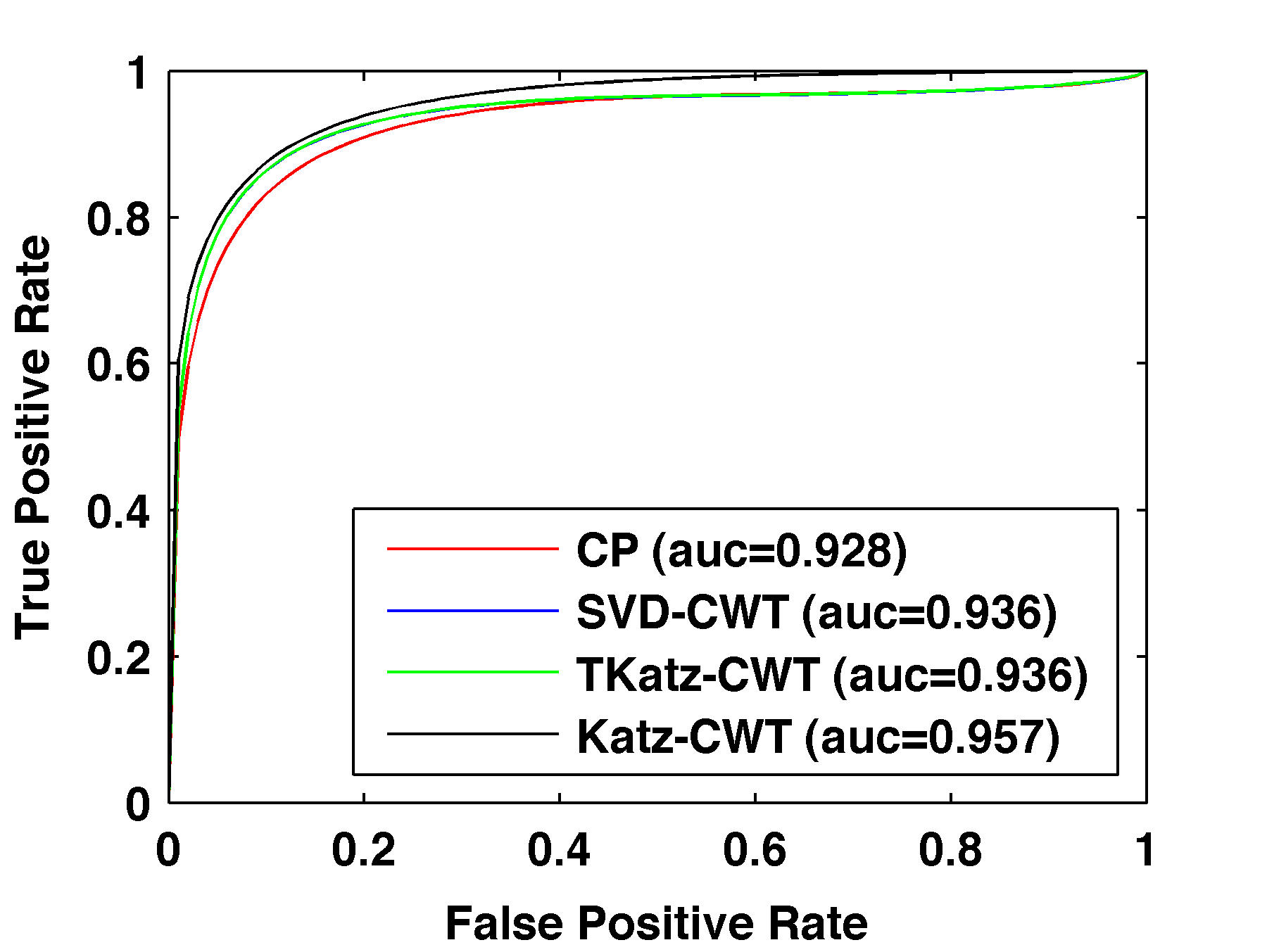}
    \label{fig:roc_all}
  }
  \subfloat[Prediction of \emph{new} links in the test sets.]
  {
    %Trim = amount to cut off on LEFT, BOTTOM, RIGHT, TOP
    \includegraphics[width=.45\textwidth,trim=0 0 0 0,clip]{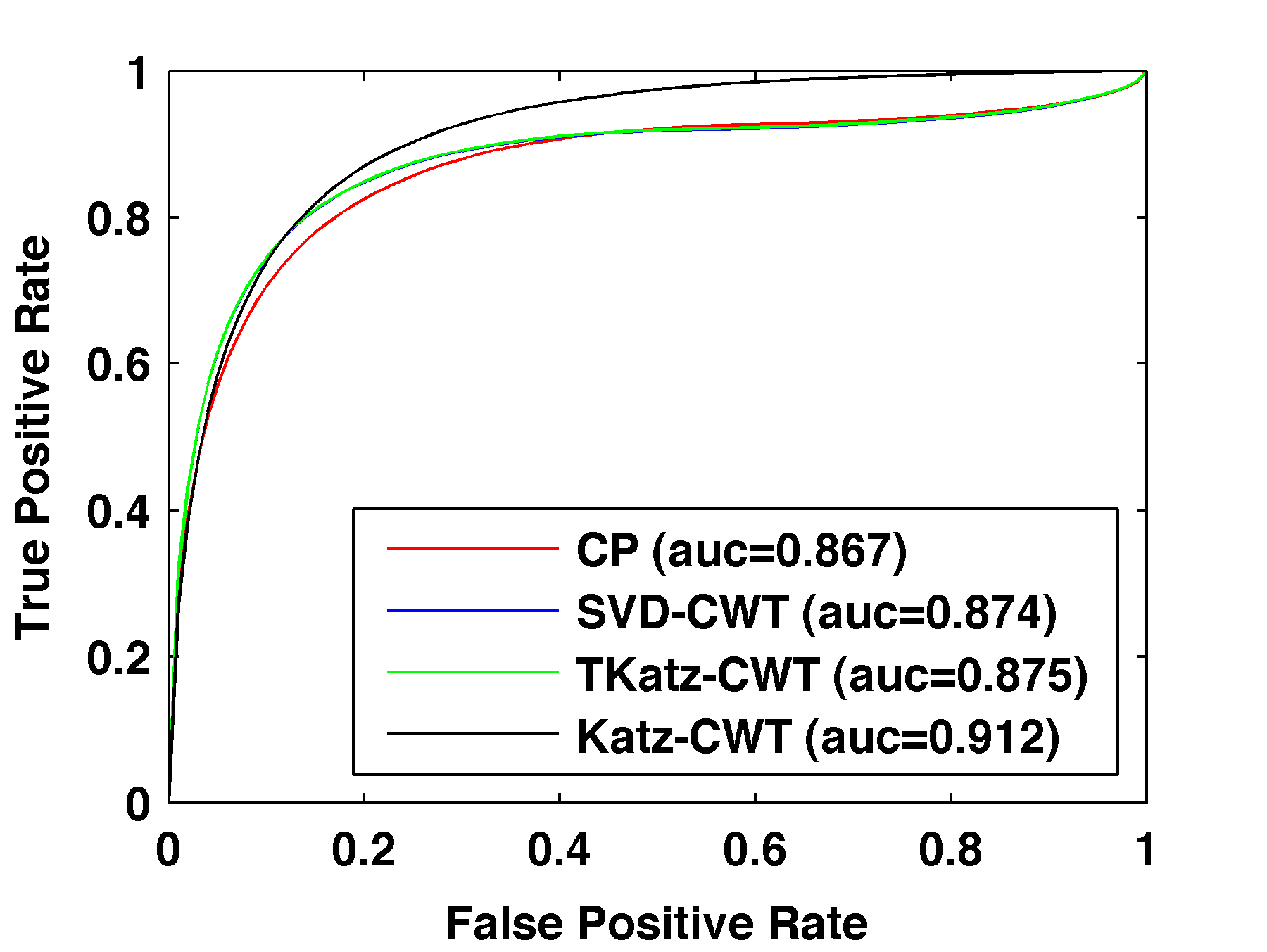}
    \label{fig:roc_unseen}
  }
  \caption{Average ROC curves showing the performance of link prediction methods across all training/test sets.}
  \label{fig:ROC}
\end{figure}

In \Fig{ROC}, we show the ROC (receiver operating characteristic) curves;
for the purposes of clarity, we omit the CT results.
When predicting all links,
\Fig{roc_all} shows that all methods perform similarly initially,
but Katz-CWT is best as the false positive rate increases.
TKatz-CWT and TSVD-CWT are only slightly worse
than Katz-CWT. Finally, CP starts having false-positives earlier than
the other methods. \Fig{roc_unseen} shows the behavior for just the
new links. In this case, the relative performance of
the algorithms is mostly unchanged.

In order to understand the behavior of different link predictors
better, we also compute how many \emph{correct links} (true positives)
are in the top 1000 scores predicted by each method. \Tab{Top1000} shows that
CP, TSVD-CWT, TKatz-CWT and Katz-CWT achieve close to 75\% accuracy
over all links.  The accuracy of the methods goes down to
10\% or less when we remove all previously seen links from the test set, but
this is still very significant.
Although 10\% accuracy may seem low, this is still two orders of
magnitude better than what would be expected by chance (0.1\%) due to high
imbalance in the data (see the last column of \Tab{Data}).
Note that the best
methods in terms of AUC, i.e., Katz-CT and Katz-CWT, perform worse
than CP, TSVD-CWT and TKatz-CWT for predicting new links.
We also observe that CP is among the
best methods when we look at the top predictions even if it starts
giving false-positives earlier than other methods.

\begin{table}[!t]
  \footnotesize
  \centering
  \setlength{\tabcolsep}{4pt}
  \subfloat{
    \footnotesize
  \begin{tabular}{| c | c c c c c c c|}
  \hline
    Test &  CP & TSVD & TSVD & TKatz & TKatz & Katz & Katz \\
    Year & & -CT & -CWT & -CT & -CWT & -CT & -CWT \\ \hline \hline
    \multicolumn{8}{|c|}{\bf All Links} \\ \hline \hline
 $2001$  &671 &617 &685 &617 &686 &625 &709\\
 $2002$  &668 &660 &674 &659 &674 &658 &716\\
 $2003$  & 723 &697 &743 &693 &745 &715 &754\\
 $2004$  & 783 &726 &777 &721 &776 &719 &774\\
 $2005$  &755 &716 &776 &720 &775 &700 &781\\
 $2006$  & 807 &729 &801 &731 &800 &698 &796\\
 $2007$  &  721 &681 &755 &687 &754 &647 &724\\  \hline
    \multicolumn{1}{|c|}{\bf Mean} & 733 &689 &744 &690 &744 &680
    &750\\ \hline
    \hline
    \multicolumn{8}{|c|}{\bf New Links} \\ \hline \hline
$2001$  &87 &80 &104 &80 &104 &51 &56\\
$2002$  &97 &84 &124 &84 &124 &74 &81\\
$2003$  &78 &80 &96 &75 &97 &55 &62\\
$2004$  &99 &79 &105 &81 &105 &57 &69\\
$2005$  &116 &89 &117 &88 &117 &58 &69\\
$2006$  &91 &77 &110 &77 &109 &63 &70\\
$2007$  &83 &71 &95 &73 &99 &43 &42\\  \hline
    \multicolumn{1}{|c|}{\bf Mean} & 93 &80 &107 &80 &107 &57 &64\\ \hline
  \end{tabular}
}
   \caption{Correct predictions in top 1000 scores.}
  \label{tab:Top1000}
\end{table}

The matrix-based methods (TSVD, TKatz and Katz) are quite fast relative to the CP method. Average timings across all training sets are
TSVD-CT/CWT and TKatz-CT/CWT: 61 sec.; Katz-CT: 80 sec.; Katz-CWT: 74 sec.; and CP: 1300 sec. Note that for all but the Katz method, the times reflect the computation of models using ranks of $K=10,20,...100$.  For the Katz method, a full SVD decomposition is computed, thus making the method too computationally expensive for very large problems. 

%%% Local Variables:
%%% mode: latex
%%% TeX-master: "paper"
%%% TeX-command-default: "PDFLaTeX"
%%% End:
     % Numerical results
\def\Train{^{\rm [train]}}
\def\Test{^{\rm [test]}}
\section{Experiments with Link Prediction for Multiple Time Steps}
\label{sec:numerical2}

As might be evident from \Figs{components}{c46}, the utility of tensor
models is in their ability to reveal patterns in time. In this section,
we use synthesized data to explore the ramifications of temporal
predictions in situations where there are several different periodic
patterns in the data. In
particular, \Fig{c46} shows an every other year pattern in the data,
but this is not exploited in the predictions in the last section based
on the heuristic score in \Eqn{score-cp}. In the DBLP data set,
periodic time profiles were few and did
not have noticeable impact on performance. Here we simulate the
type of data where bringing the time pattern information into play is
crucial
for predictions. Our goal is to use the periodic information,
predicting even further out in
time. For example, if we assume that our training data is for times
$t=1,\dots,T$ and that the period in our data is of length $L$, then
we can predict connections for time periods $T+1$
through $T+L$. In this section, we present results of experiments
involving link prediction for multiple time steps. All experiments were
performed using Matlab 7.9 on a Linux  Workstation (RedHat 5.2)
with 2 Quad-Core Intel Xeon 3.0GHz processors and 32GB RAM.

\subsection{Data}

We generate simulated data that shows connections between two sets
of
entities over time. In the DBLP data, the entity sets were authors and
conferences and each time $t=1,\dots,T$ corresponded to one year of
data. In our simulated example, we assume that each time $t$
corresponds to one day and that the temporal profiles correspond
roughly
to a seven-day period ($L=7$). We might assume that the entities are
users and services. For example, most major service providers
(Yahoo,
Google, MSN, etc.) have a front page. We may wish to predict which
users will likely be connecting to which services over time. For
example, it may be that there are groups of people that check the
National, World, and Business News on weekdays; another group that
checks Entertainment Listings on weekends; a group that checks
Sports
Scores on Mondays; another group uses the service for email every
day,
etc. Each of these groupings of users and services along with
corresponding temporal profile can be represented by a single
component in the tensor model. The motivation for link prediction are
many-fold. We might want to characterize the temporal patterns so that
we know how often to update the services or when to schedule down
time. We may use prediction to cache certain data,
to better direct advertisements to users, etc.  In addition to the
business model/motivation, we could also consider the application of
cybersecurity using analysis of network traffic data. In this case,
the goal is to determine which computers are most likely to contact
which other computers. Predicted links could be used to find
anomalous
or malicious behavior, proactively load balance, etc. In all of these
applications, accurately predicting links over multiple time steps is
crucial.

We assume that our data can be modeled by $K=10$ components and
generate our training and testing tensors as follows.
\begin{enumerate}
\item
Matrices
$\M{A}$ and $\M{B}$ of size $M \times K$ and $N \times K$ are
``entity
participation'' matrices. In other words, column $\MC{A}{k}$
(resp. $\MC{B}{k}$) is the
vector of participation levels of all the entities in component $k$.
In our tests, we use {$M=500$ and
  $N=400$.} Each row is generated by choosing between 1
and $K$ components for the entity to participate in where the
probability of participating in at least $k+1$ components is $1 - 4^{k}$.
An example of the distribution of participation is shown in \Fig{cdf};
note that most
entities participate in just one component. Once the number of
components
is decided, the specific components that a given entity participates in
are chosen uniformly at random. The strength of participation is
picked uniformly at random between 1 and 10. Finally, the columns of
$\M{A}$ and $\M{B}$ are normalized to length 1.

\begin{figure}[t!]
  \centering
  %Trim = amount to cut off on LEFT, BOTTOM, RIGHT, TOP
  \includegraphics[trim=0 5 0 10,clip,width=3in]{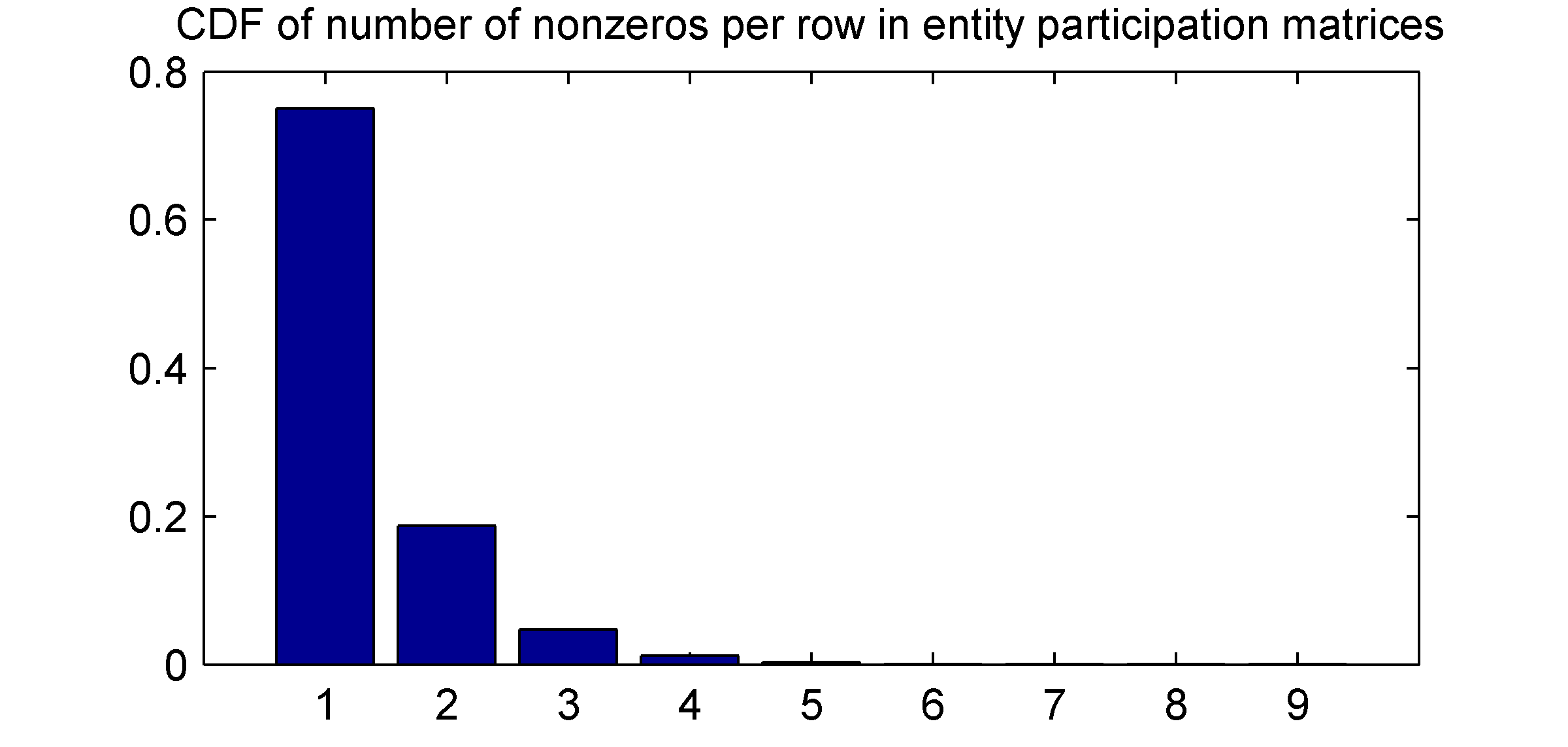}
  \caption{Example distribution of the number of
    nonzeros per row in entity participation matrices.}
  \label{fig:cdf}
\end{figure}

\item In the third mode, the matrix corresponds to time. We assume
  we have $P$ periods of training data, so that we have $T=LP$
time observations
  to train on. In our tests, we use $L=7$ and $P=10$.
  We also assume that we have 1 period of testing
  data. Therefore, we generate a matrix $\M{C}$ of size $L(P+1) \times
  K$, which will be divided into submatrices $\M{C}\Train$ and
  $\M{C}\Test$. Each column of $\M{C}$ is temporal data of length
  $L(P+1)$ with a repeating period of length $L=7$. We use the
  periodic patterns shown in
  \Fig{patterns}. For example, the first patten corresponds to a
  weekday pattern, and the seventh pattern corresponds to
  Tuesday/Thursday activities. Patterns 1 and 5 are the
  same but correspond to different sets of entities and so will still
  be computable in our experiments.

  \begin{figure}[!t]
    \centering
    %% Trim = amount to cut off on LEFT, BOTTOM, RIGHT, TOP
    \includegraphics[trim=0 20 0 20,clip,width=4in]{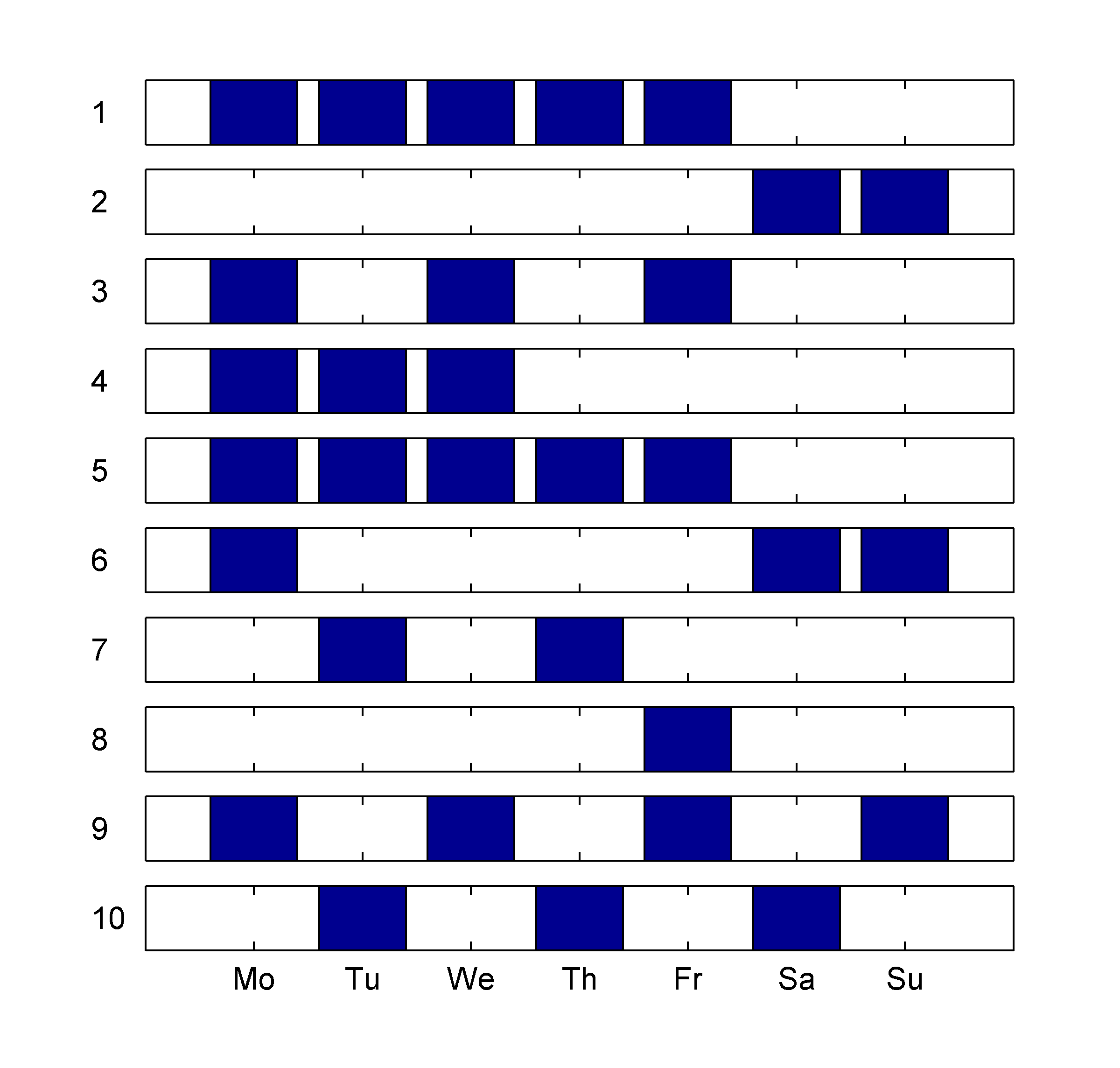}
    \caption{Weekly patterns of users in simulated data.}
    \label{fig:patterns}
  \end{figure}

  The creation of the temporal data is shown in \Fig{pattern_example}.
  As shown at the top, in order to generate the temporal data of
  length $L(P+1)$, we repeat the temporal pattern $P+1$ times. Next,
  for each of the $K=10$ components, we randomly adjust the data to
be
  increasing, decreasing, or neutral. On the left in the middle plot,
  we show a decreasing pattern and on the right an increasing pattern.
  Finally, we add 10\% noise, as shown in the lower plots.
  The final matrix $\M{C}$ is column normalized. The first $T=LP=70$
rows
  become $\M{C}\Train$ and the last $L=7$ rows become $\M{C}\Test
$.

  \begin{figure}[htbp]
    \centering
    %% Export set-up:
    %% Size 4" x 4.16", Check "expand axis" box
    %% Rendering: 600 dpi
    %% Fonts: Custom size fixed at 10 points
    \subfloat[Pattern 3 and a decreasing trend]%
    {\includegraphics[trim=25 20 20 0,clip,height=2.5in]
      {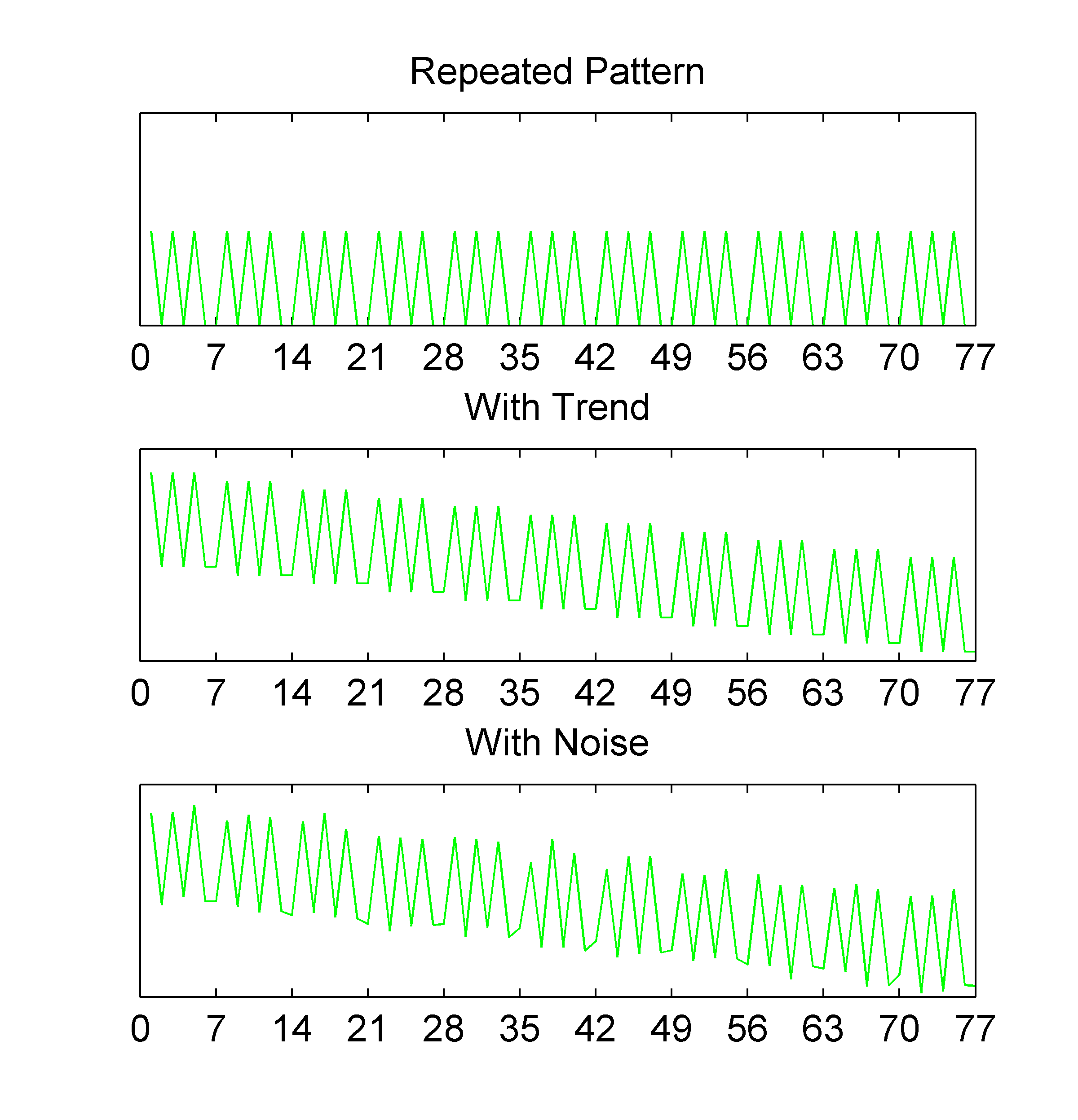}}
    ~~
    \subfloat[Pattern 1 and an increasing trend]%
    {\includegraphics[trim=25 20 20 0,clip,height=2.5in]
      {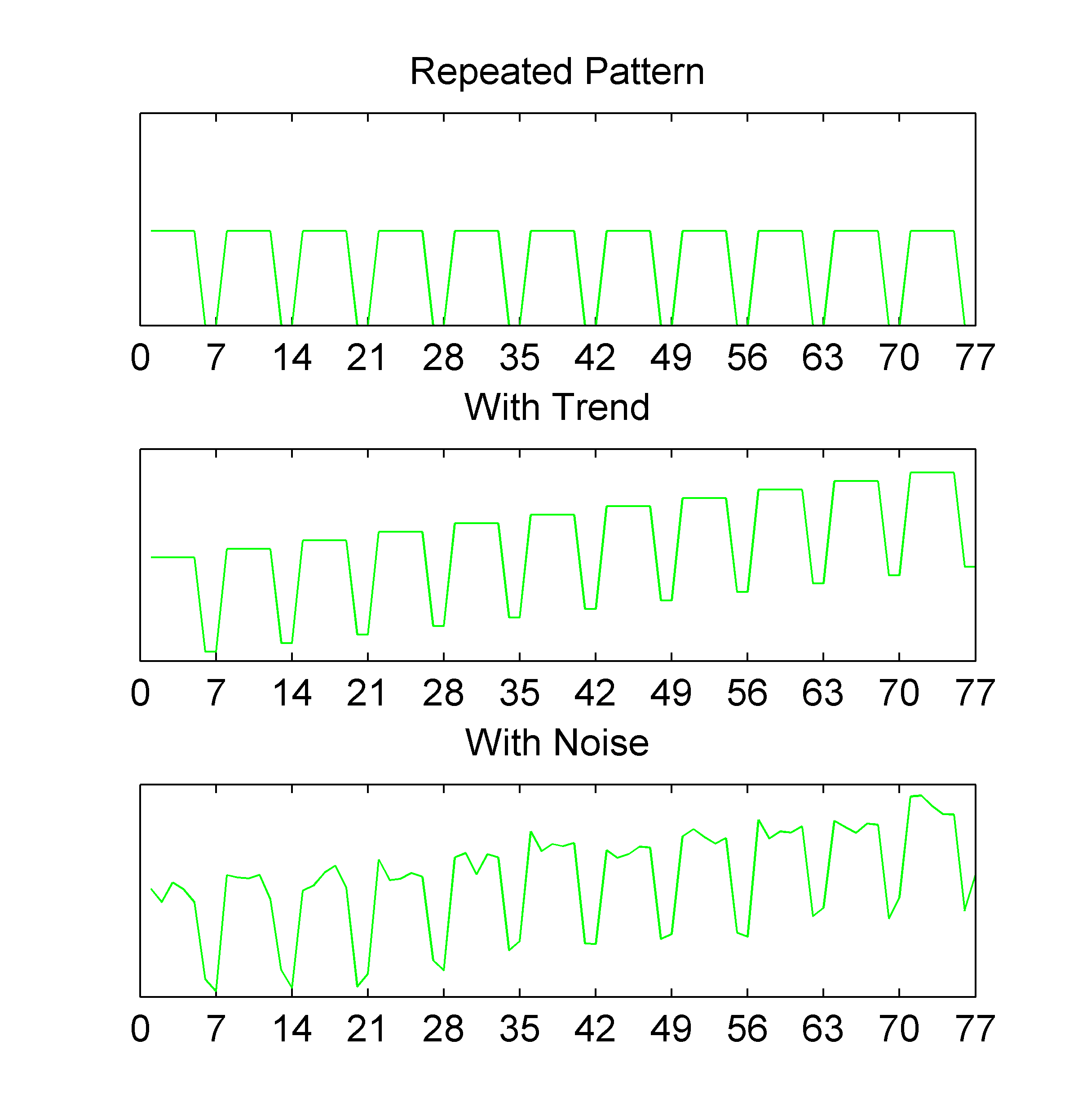}}
    \caption{Steps to creating the temporal data.}
    \label{fig:pattern_example}
  \end{figure}

\item We create noise-free versions of the training and testing tensors
by
  computing
  \begin{displaymath}
    \T{Z}\Train = \sum_{k=1}^K \MC{A}{k} \Oprod \MC{B}{k} \Oprod
    \MC{C}{k}\Train
    \quad\text{and}\quad
    \T{Z}\Test = \sum_{k=1}^K \MC{A}{k} \Oprod \MC{B}{k} \Oprod
\MC{C}{k}\Test.
  \end{displaymath}
  In order to make the problem challenging, we significantly degrade
  the training data in the next two steps.
\item We take $p_{\rm swap}=50\%$ of the $p_{\rm top}=25\%$ largest entries in $\T{Z}
\Train$ and randomly swap them with other entries in the
tensor
  selected uniformly at random. This has the effect of removing some
  important data (i.e., large entries) and also adding some spurious
  data.
\item Finally, we add $p_{\rm rand}=10\%$ percent standard normal
noise to every entry of
  $\T{Z}\Train$.
\end{enumerate}

For each problem instance, we compute a CP decomposition of $\T{Z}
\Train$ which has had
large entries swapped and noise added. Then we use the resulting
factorization to
predict the largest entries of $\T{Z}\Test$.

\subsection{Methods and Parameter Selection}
\label{sec:setup}
The goal of this study is to predict the significant entries in
$\T{Z}\Test$. Therefore, without loss of generality, we treat
all nonzeros in the test tensor as ones (i.e., positive
links) and the rest as zeros (i.e., no link). This results in 15\% positives.

We consider the CP model (with $K=10$ components) using the
forecasting scoring model described in \Sec{forecasting}. The model parameters
for level, trend and seasonality are set to $0.2$ in the Holt-Winters method. We compute
the CP model via the CPOPT approach as described in
\cite{SAND2009-0857}. This is an optimization approach using the
nonlinear conjugate gradient method. We set the stopping tolerance on
the normalized gradient to be $10^{-8}$, the maximum number of
iterations to 1000, and the maximum number of function values to
10,000. The models could have been computed using ALS as in the
previous section, but here we use a different technique for variety.

The matrix
methods described in \Sec{matrix} are not appropriate for this data
because simply collapsing the data will obviously not work. Moreover,
there is no clear methodology for predicting out in time. The best we
could possibly do would be to construct a matrix model for each day in
the week (i.e., one for each $\ell = 1,\dots,L$).  Such an approach,
however, is not parsimonious and would quickly become prohibitive as the number of models
grew. Moreover, it would be extremely difficult to assimilate results
across the different models for each day.

Therefore, as a comparison, we use the largest values from the most
recent period. In MATLAB notation, the predictions are based on
$\T{Z}\Train(:,:,L(P-1)+1:LP)$, i.e., the last $L$ frontal slices of
$\T{Z}\Train$.  The highest values in that period are predicted to
reappear in the testing period.
We call this the \emph{Last Period} method.
Under the scenario considered here, this is an
extremely predictive model when the noise levels are low. As we
randomly swap data (reflective of random additions and deletions as would be expected in any real world data set),
however, the performance of this model degrades.

\subsection{Interpretation of CP Results and Temporal Prediction}

As mentioned above, we compute the CP factorization of a noisy
version
of $\T{Z}\Train$ of size $500 \times 400 \times 70$. Because the data
is noisy (swapping 50\% of the 25\% largest
entries and adding 10\% Gaussian noise), the fit of the model to the
data is not perfect. In fact,
the percentage of the $\T{Z}\Train$ that is described by the model%
\footnote{The percentage of the data described by the model is
calculated as $1 - \| \T{M} - \T{Z}\Train \| / \|\T{Z}\Train \|$ where
$\T{M}$ is the CP model.}
is only about 50\% (averaged over 10 instances).

However, the underlying low-rank structure of the data gives us hope
of recovery even in the presence of excessive errors. We can see this
in terms of the \emph{factor match score} (FMS), which is defined as
follows.
Let the correct and computed factorizations be given by
\begin{displaymath}
  \sum_{k=1}^K \lambda_k \; \MC{A}{k} \Oprod \MC{B}{k}
  \Oprod \MC{C}{k}
  \quad\text{and}\quad
  \sum_{k=1}^{K} \bar\lambda_k \; \MbarC{A}{k} \Oprod \MbarC{B}{k}
  \Oprod \MbarC{C}{k},
\end{displaymath}
respectively. Without loss of generality, we assume that all the
vectors have been scaled to unit length and that the
scalars are positive. Recall that there is a permutation
ambiguity, so all possible matchings of
components between the two solutions must be considered.
Under these conditions, the FMS is defined as
\begin{equation}
\label{eq:FMS}
  \text{FMS} = \max_{\sigma\in\Pi(K)}
  \frac{1}{K} \sum_{k=1}^K
  \left(
    1 - \frac{|\lambda_r -  \bar\lambda_{\sigma(k)}|}
    {\max\{\lambda_r,\bar\lambda_{\sigma(k)}\}}
  \right)
  |\MC{A}{k}\Tra \MbarC{A}{\sigma(k)}|
  |\MC{B}{k}\Tra \MbarC{B}{\sigma(k)}|
  |\MC{C}{k}\Tra \MbarC{C}{\sigma(k)}| .
\end{equation}
The set $\Pi(K)$ consists of all permutations of 1 to $K$. In our case,
we just
use a greedy method to determine an appropriate permutation. The
FMS can be
between 0 and 1, and the best possible FMS is 1. For the problems mentioned above, the FMS
scores are around 0.52 on average; yet the predictive power of the CP model is very good
(see \Sec{periodicResults}).

\begin{figure}[!t]
  \begin{center}
    %% Trim = amount to cut off on LEFT, BOTTOM, RIGHT, TOP
    \includegraphics[trim=50 60 50 30,clip,width=\textwidth]%
    {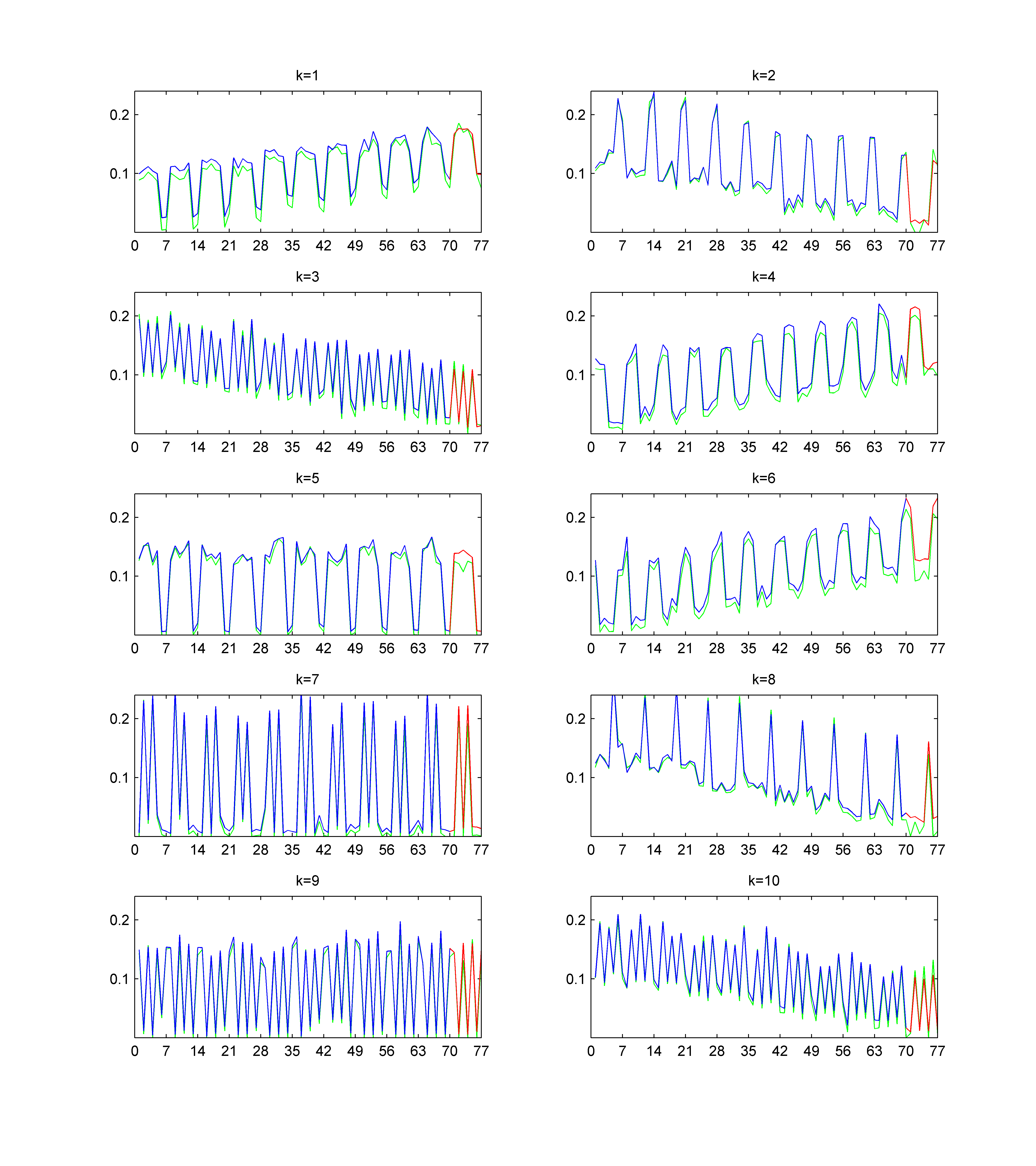}
  \end{center}
  \caption{Temporal patterns using Holt-Winters forecasting. The green
    line is the ``true'' data. The blue line is the temporal pattern
    that is computed by CP. The red line is the pattern that is
    predicted by Holt-Winters using the temporal pattern computed by
CP.}
\label{fig:temporal_trend}
\end{figure}

\Fig{temporal_trend} shows the ten different temporal patterns in the
data per the $\MC{C}{k}$ vectors. The
green line is the original unaltered data, the first portion of
which is used to generate the training data. The blue line is the
pattern that is extracted via the CP model. Note that it is generally
a good fit to the true data shown in green. Finally, the red line is
the prediction generated by the Holt-Winters method using the pattern
extracted by CP (i.e., the blue line). The red data corresponds to the
$\bm \gamma_k$ vectors used in \Eqn{score2-cp} for link
prediction. The results of the link prediction task are discussed
in the next subsection.

\subsection{Link Prediction Results}
\label{sec:periodicResults}

\begin{figure}[!t]
\begin{center}
\includegraphics[width=.75\textwidth]{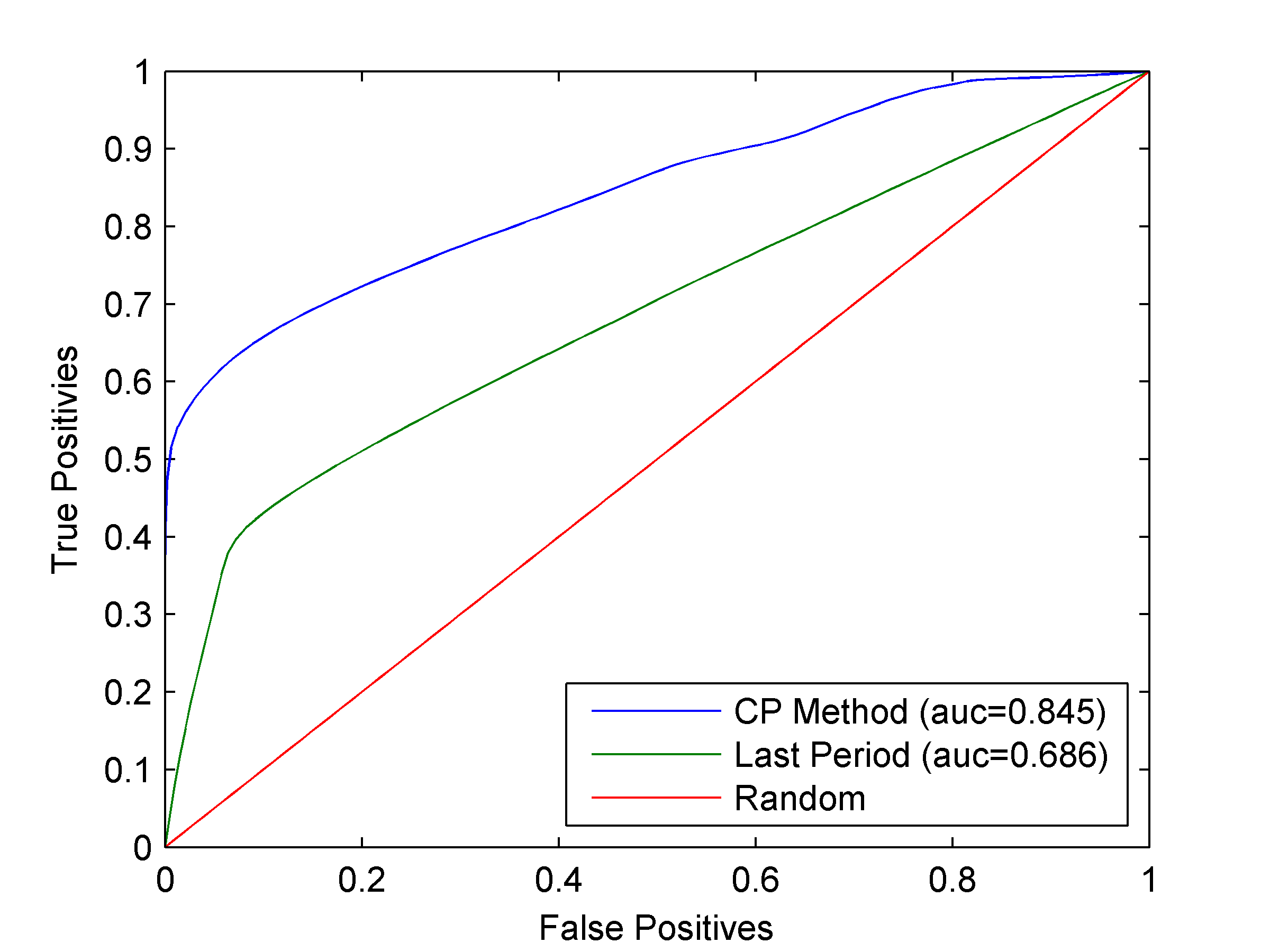}
\end{center}
\caption{ROC curves and AUC scores for typical problem of predicting links seven steps forward in time.}
\label{fig:auc_trend}
\end{figure}

In \Fig{auc_trend}, we present typical ROC curves for the predictions of the next
$L=7$ time steps (one period) for a problem instance generated using the
procedure mentioned above. As described in \Sec{setup}, our goal
is to predict the nonzero entries in the testing data
$\T{Z}\Test$ based on the CP model and the score in
\Eqn{score2-cp}. We compare with the predictions based on the last
period in the data. Despite the high level of noise, the CP method is
able to get an AUC score of $0.845$, which is much
better than then the ``Last Period'' method's score of $0.686$.
We also considered the accuracy in the first 1,000 values
returned. The CP-based method is 100\% accurate in its top 1,000
scores whereas the ``Last Period'' method is only 70\% accurate.

To investigate the effect of noise on the performance of the CP and
Last Period methods, we ran several experiments varying the different
amounts of noise ($p_{\rm top}$, $p_{\rm swap}$, and $p_{\rm rand}$). For each
experiments, we fixed all but one type of noise, varying the remaining
type of noise. The fixed values for each type of noise were the same as
those for the experiment described above,
while varying $p_{\rm top}$ up to 30\%,
$p_{\rm swap}$ up to 80\%, and $p_{\rm rand}$ up to 40\%.
For each level of noise, we generated 10 instances of
$\T{Z}\Train$ and $\T{Z}\Test$ and computed the average AUC and
percentage of links correctly predicted in the top 1000 scores.

\begin{figure}[!tp]
    \centering
    \subfloat[AUC]%
    {\includegraphics[width=.4\textwidth]
      {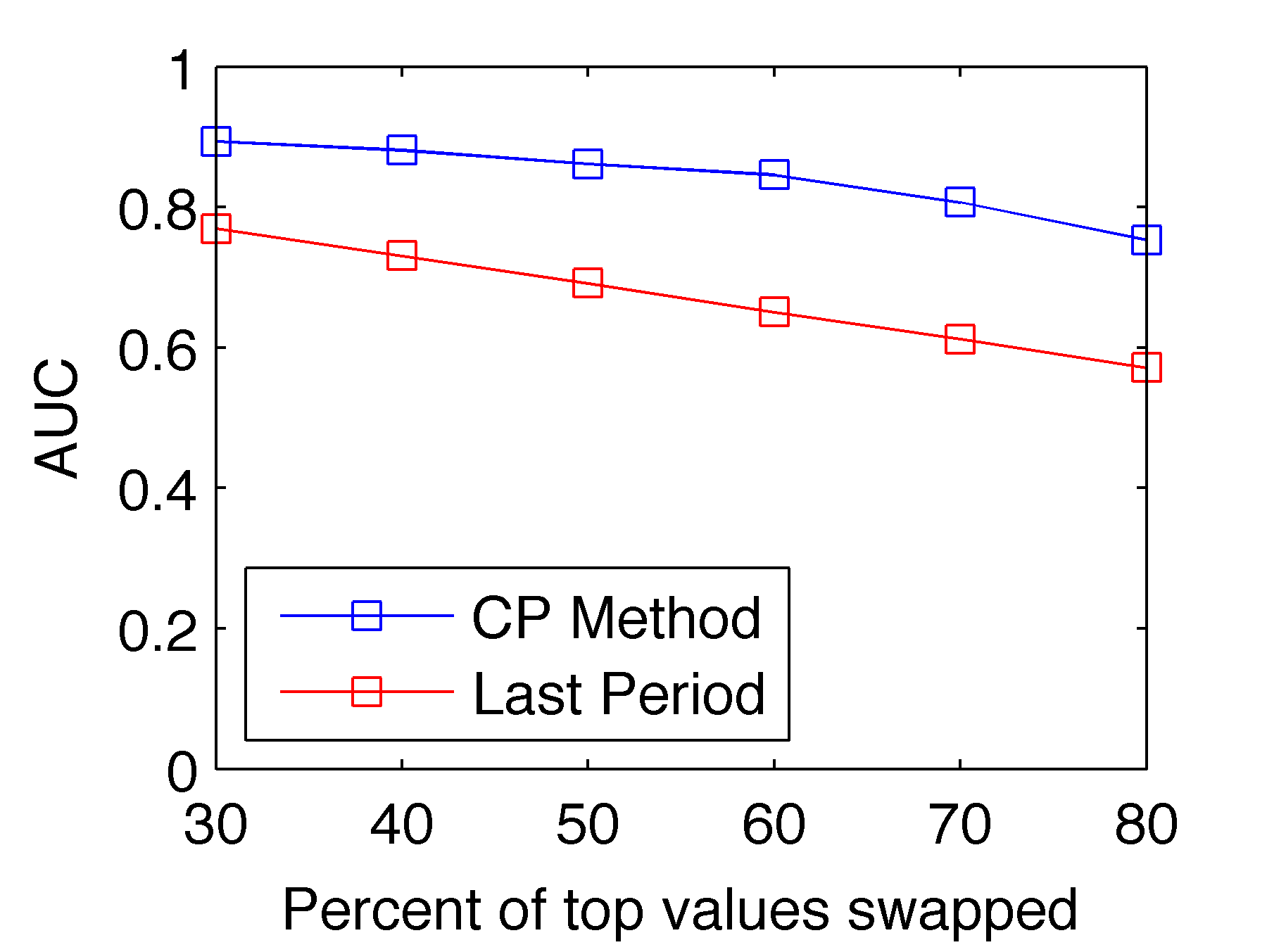}}
    ~~
    \subfloat[Percent Correct in Top 1000]%
    {\includegraphics[width=.4\textwidth]
      {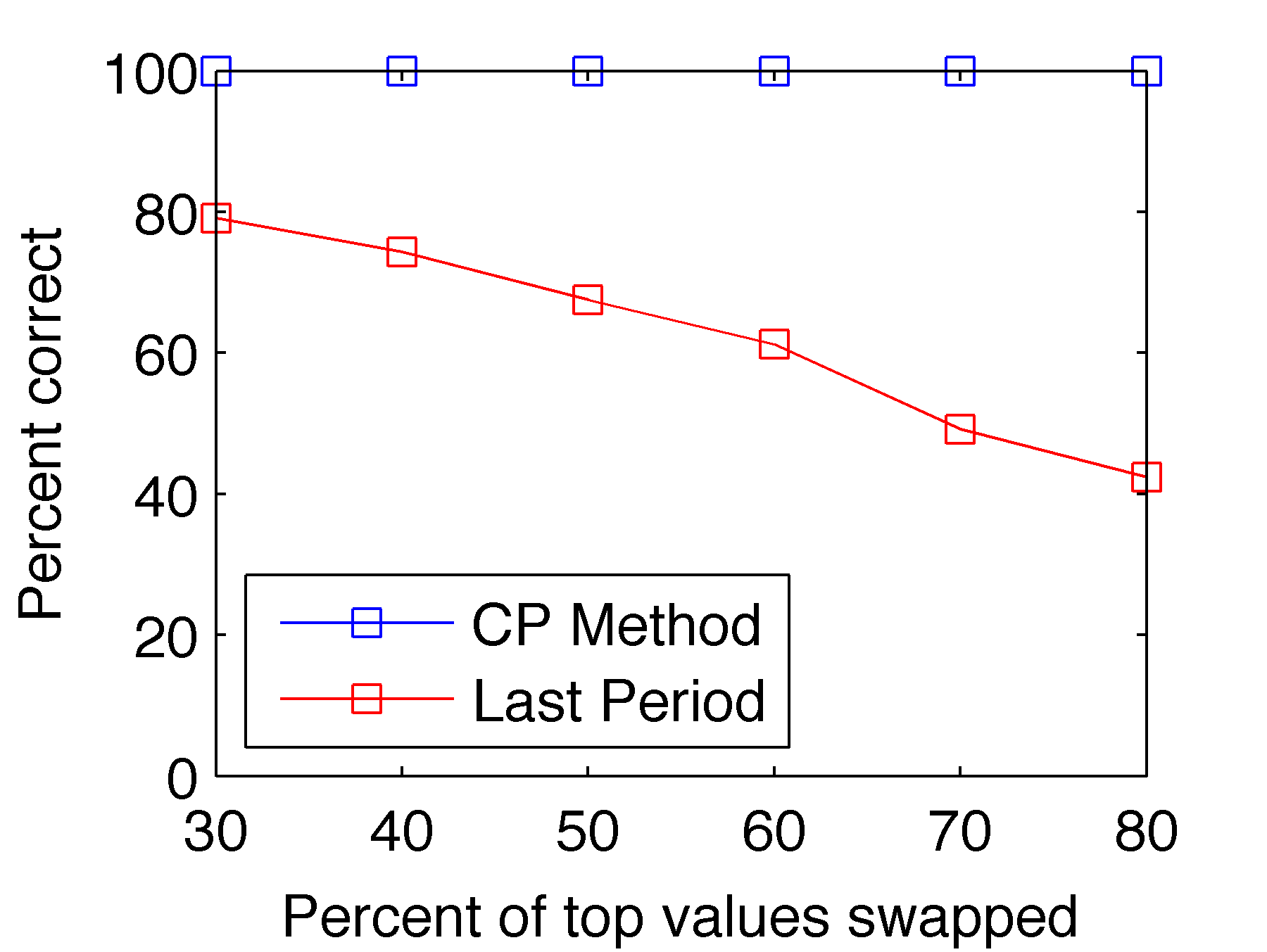}}
    \caption{Impact of varying swapping noise ($p_{\rm swap}$)
averaged over 10 runs per noise value.}
    \label{fig:swapNoise}
  \end{figure}

\begin{figure}[!tp]
    \centering
    \subfloat[AUC]%
    {\includegraphics[width=.4\textwidth]
      {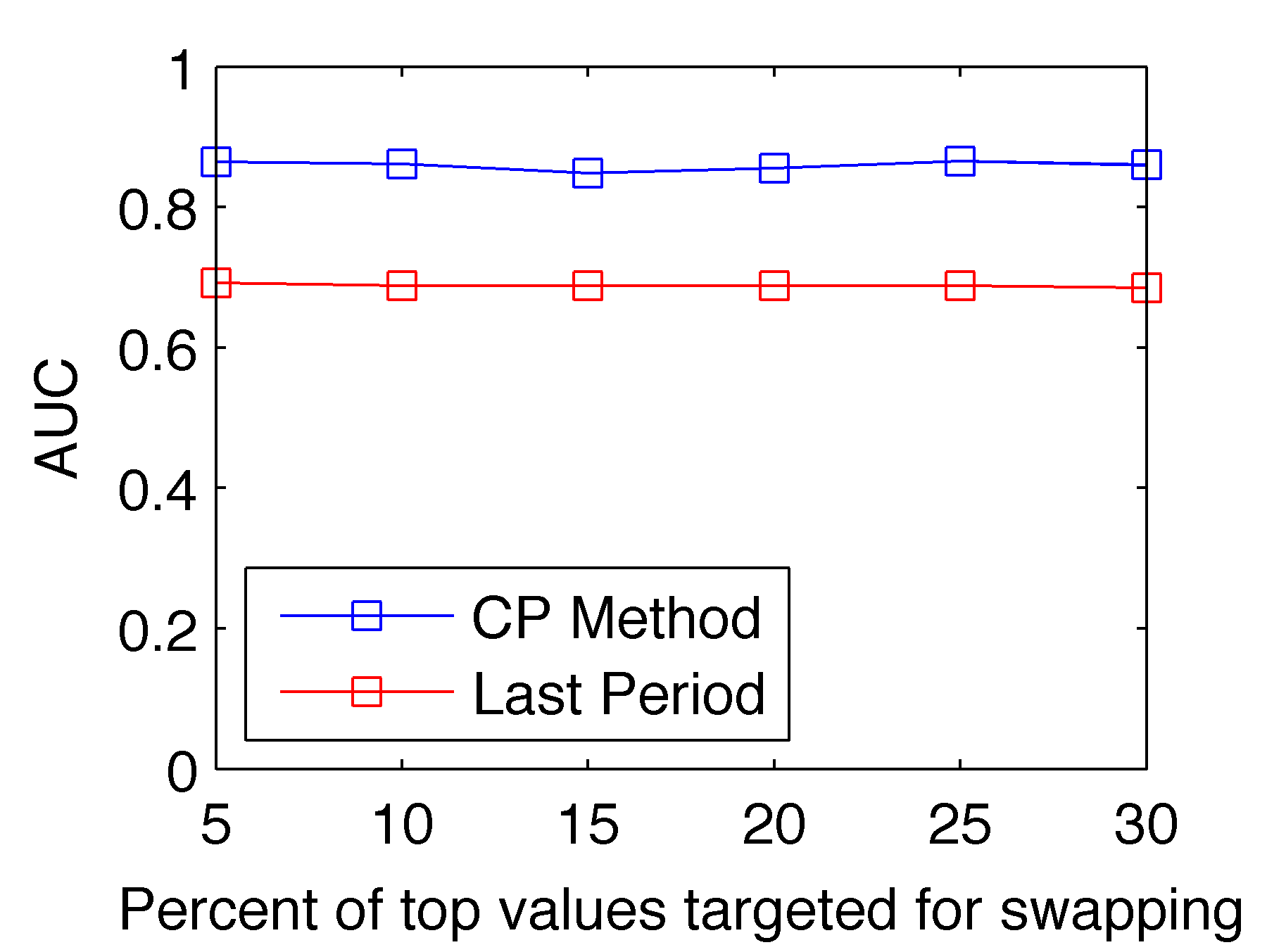}}
    ~~
    \subfloat[Percent Correct in Top 1000]%
    {\includegraphics[width=.4\textwidth]
      {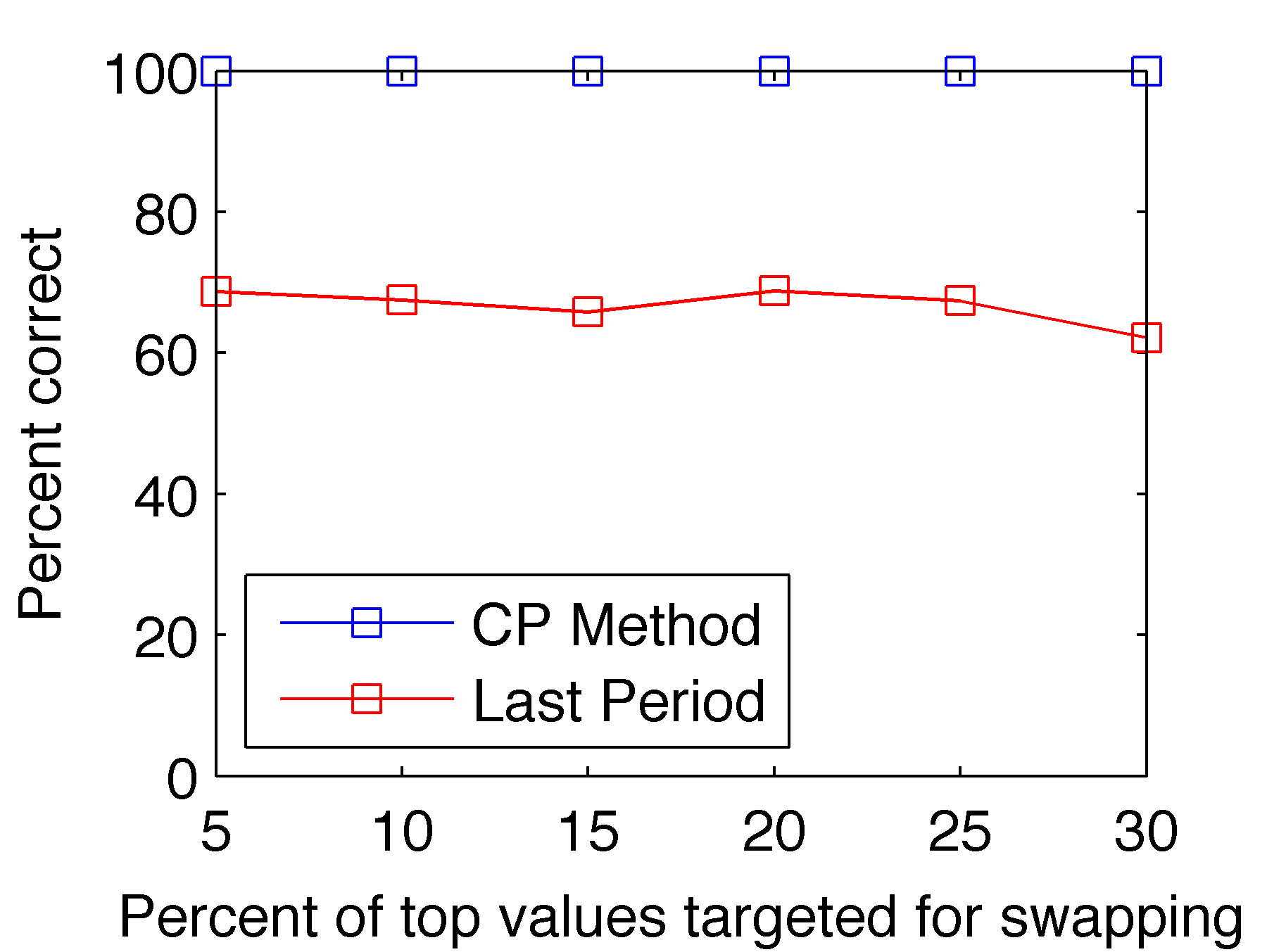}}
    \caption{Impact of varying random noise ($p_{\rm top}$)
averaged over 10 runs per noise value.}
    \label{fig:topValThresh}
  \end{figure}

\begin{figure}[!tp]
    \centering
    \subfloat[AUC]%
    {\includegraphics[width=.4\textwidth]
      {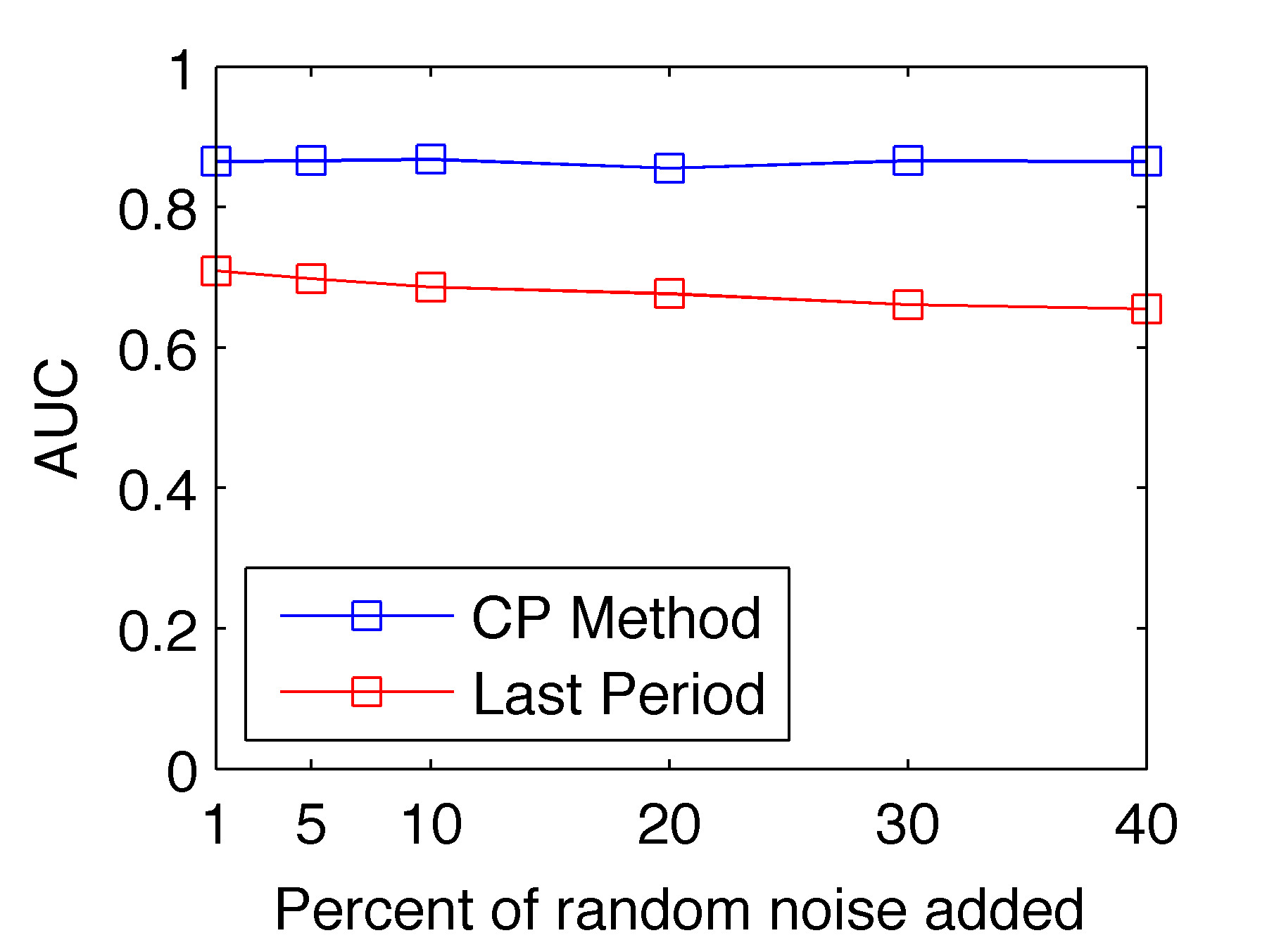}}
    ~~
    \subfloat[Percent Correct in Top 1000]%
    {\includegraphics[width=.4\textwidth]
      {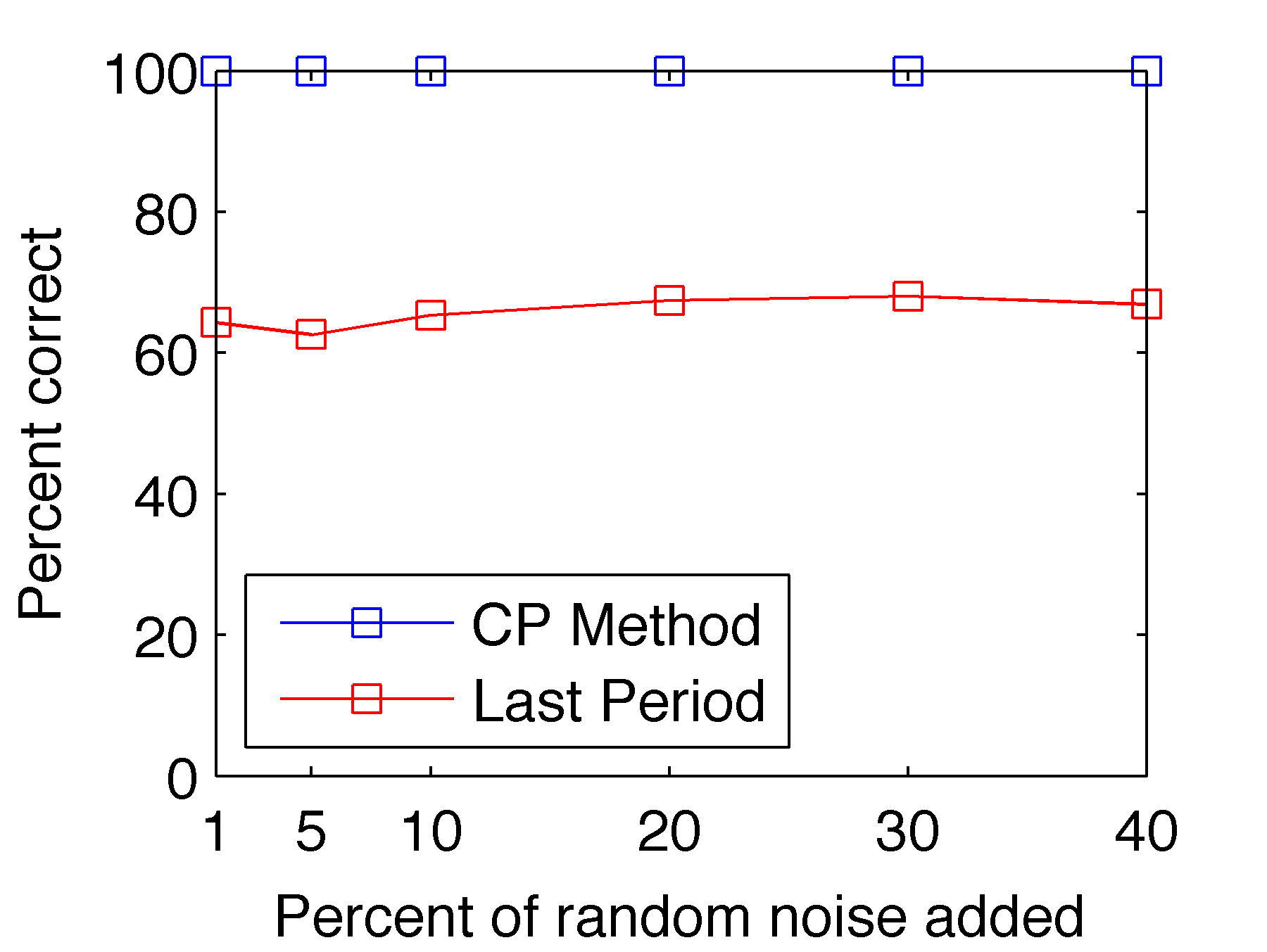}}
    \caption{Impact of varying random noise ($p_{\rm rand}$)
averaged over 10 runs per noise value.}
    \label{fig:randomNoise}
  \end{figure}

\Fig{swapNoise} shows the results of the experiments when varying
$p_{\rm swap}$. As expected, AUC values decrease as the number of the most significant
links in the training data being swapped randomly is increased. However,
there is a clear advantage of the CP method over the Last Period method as depicted in \Fig{swapNoise}a. Note also that in \Fig{swapNoise}b, we see that even as the number of randomly swapped significant links is increased in the training data, the top 1000 scores predicted with the CP method are all correctly identified as links. For the Last Period method, performance decreases as $p_{\rm swap}$ increases, highlighting a clear advantage of the CP method. \Fig{topValThresh} and \Fig{randomNoise} present the results for the experiments where $p_{\rm top}$ and $p_{\rm rand}$ were varied, respectively. In both sets of experiments, the CP method consistently performed better than the Last Period method in terms of both AUC and correct predictions in the top 1000 scores for each method. However, no significant changes were detected across the different levels of $p_{\rm top}$ and $p_{\rm rand}$.

The key conclusions from these experiments is that the CP model performs extremely well even when a large percentage of the strongest signals (i.e., link information) in the training data is altered. Such robustness is crucial for applications where noise or missing data is common, e.g., analysis of computer network traffic where data is lost or even hidden in the context of malicious behavior or the analysis of user-service relationships where both user and service profiles are changing over time.

\begin{figure}[!t]
    \centering
    \subfloat[Timings]%
{\includegraphics[width=.4\textwidth]{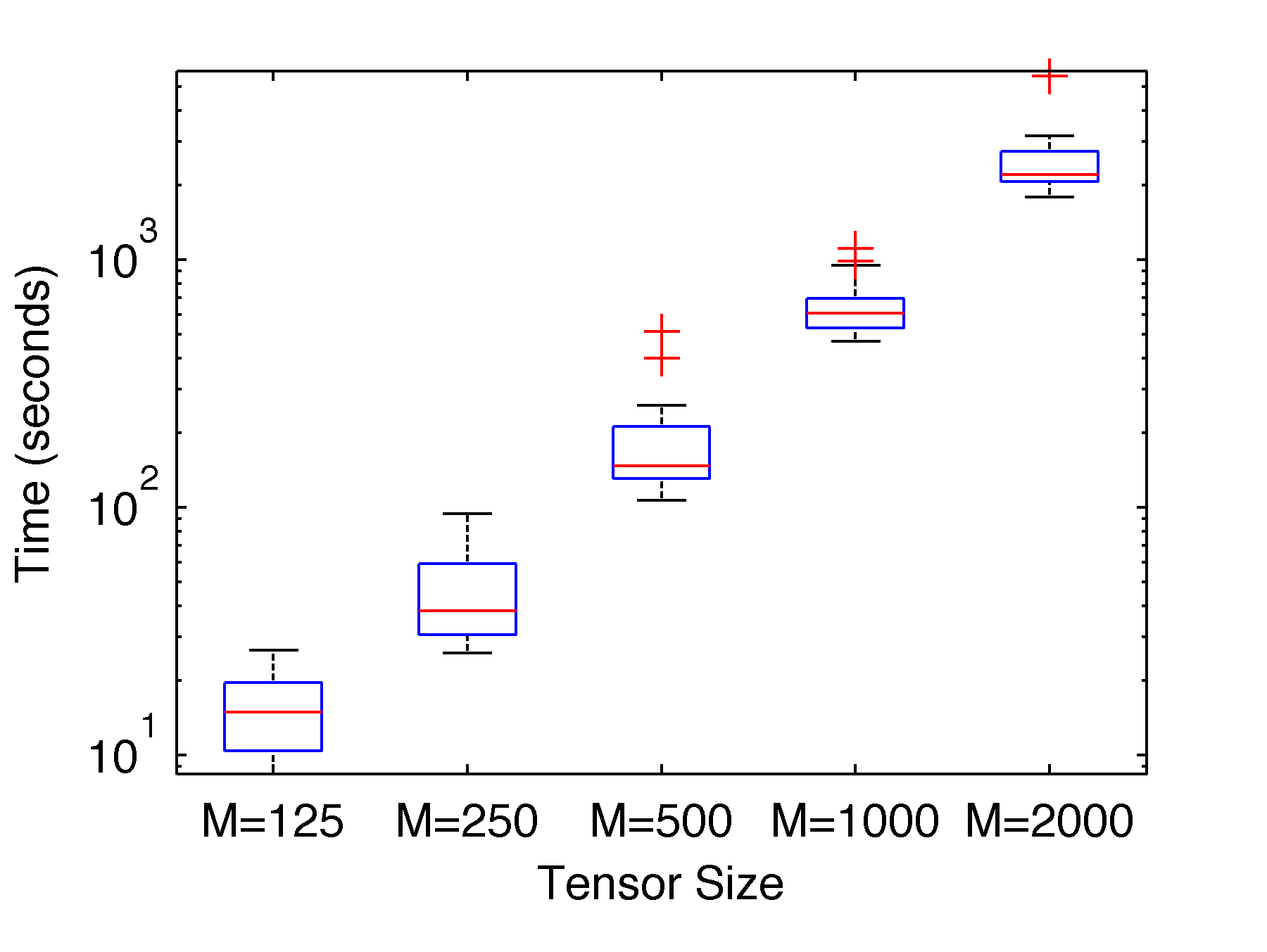}}
    ~~
    \subfloat[AUC]%
{\includegraphics[width=.4\textwidth]{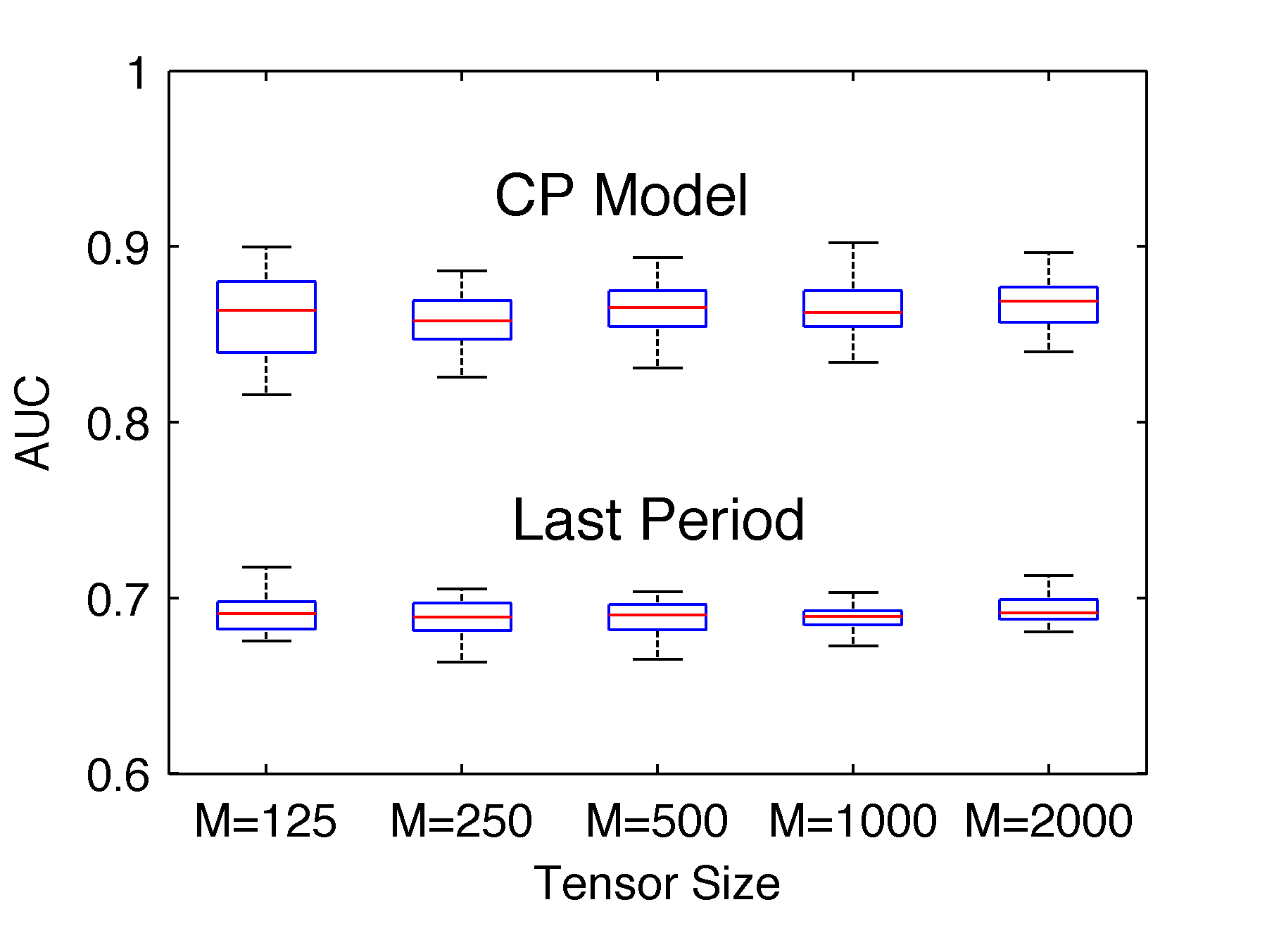}}
\caption{Results of computing 10 CP models and Holt-Winters forecasts for tensors of different sizes: $125\times100\times77$ (M=125); $250\times200\times77$ (M=250); $500\times400\times77$ (M=500); $1000\times800\times77$ (M=1000); $2000\times1600\times77$ (M=2000). Computation wall clock times in seconds are shown in (a) and AUC scores are shown in (b).}
\label{fig:timings_scalability}
\end{figure}

Results for different sizes of tensors illustrate that these conclusions hold for larger data sets as well. \Fig{timings_scalability} presents plots of (a) the time required to compute the CP models and Holt-Winters forecasts and (b) the AUC scores for predicting seven steps out in time for tensors of sizes $125\times100\times77$, $250\times200\times77$, $500\times400\times77$ (M=500), $1000\times800\times77$ (M=1000), and $2000\times1600\times77$. For each size, 10 tensors were generated using the procedures in \Sec{setup}, and the box and whiskers plots in \Fig{timings_scalability} present the median (red center mark of boxes), middle quartile (top and bottom box edges), and outlier (red plus marks) summary statistics across the experiments. These results support the conclusions above: predictions using the CP model are more accurate than those computed using the last period to predict an entire period of links.

%%% Local Variables:
%%% mode: latex
%%% TeX-master: "paper"
%%% TeX-command-default: "PDFLaTeX"
%%% LaTeX-command-style: (("." "pdflatex"))
%%% End:
     % Numerical results
\section{Related Work}
\label{sec:related}

Getoor and Diehl \citeyear{GeDi05} present a survey of link mining tasks,
including node classification, group detection, and numerous other
tasks including link prediction. Sharan and Neville \citeyear{ShNe08}
consider the goal of node classification for temporal-relational
data, suggesting the idea of a ``summary graph'' of weighted
snapshots in time which we have incorporated into this work.

The seminal work of Liben-Nowell and Kleinberg \citeyear{LiKl07} examines
numerous methods for link prediction on co-authorship networks in
arXiv bibliometric data.  However, temporal information was unused
(e.g., as in \cite{ShNe08}) except for splitting the data.
The proportion of new links ranged from 0.1--0.5\% and is thus comparable to what
we see in our data (0.05--0.07\%). According to Liben-Nowell and Kleinberg, Katz and its
variants are among the best link predictors; this observation has
been supported by other work as well \cite{Huang2005,Wang2007}.  We
note that Wang, Satuluri and Parthasarathy \citeyear{Wang2007}
use the truncated sum approximate Katz measure discussed in \Sec{katz} and
 recommend AUC as one evaluation measure for link prediction because
it does not require any arbitrary cut-off.
Rattigan and Jensen \citeyear{RaJe05} contend that the link mining problem
is too difficult, in part because the proportion of actual links is
very small compared to the number of possible links; specifically, they consider
co-author relationships in DBLP data and observe that the proportion
of new links is less than 0.01\%. (Although we also use DBLP data, we
consider author-conference links which has 0.05\% or more new links.)

Another way to approach link prediction is to treat it as a
straightforward classification problem by computing features for
possible links and using a state-of-the-art classification engine like
support vector machines. Al Hasan et al.\@ \citeyear{Mohammad2006} use this approach in the
task of author-author link prediction. They randomly pick equal sized
sets of linked and unlinked pairs of authors. They compute features
such as keyword similarity, neighbor similarity, shortest path,
etc. However, it would likely be computationally intractable to use such a
method for computing \emph{all} possible links due to the size of the
problem and imbalance between linked and unlinked pairs of authors.
Clauset, Moore, and Newman \citeyear{ClMoNe08} predict links (or
anomalies) using Monte-Carlo sampling on all possible dendrogram models of a graph.
Smola and Kondor \citeyear{SmKo03} identify connections between link prediction methods and
diffusion kernels on graphs but provide no numerical experiments to support this.

Modeling the time evolution of graphs has been considered, e.g., by
Sakar et al.\@ \citeyear{Purnamrita2007} who create time-evolving
co-occurrence models that map entities into an evolving latent space.
Tong et al.\@ \citeyear{Tong2008} also compute centrality measures on time
evolving bipartite graphs by aggregating adjacency matrices over time
in similar approaches to those in \Sec{collapse}.

Link prediction is also related to the task of collaborative
filtering. In the Netflix contest, for example, Bell and Koren
\citeyear{Bell2007} consider the ``binary view of the data'' as
important as the ratings themselves. In other words, it is important
to first predict who is likely to rate what before focusing on the
ratings. This was a specific task in KDD Cup 2007 \cite{Liu2007}.
More recent models by Koren \cite{Ko09}
explicitly account for changes in user preferences over time. And
Xiong et al.\@ \citeyear{XiChHuScCa10} propose a probabilistic tensor
factorization to address the problem of collaborative filtering over time.

Tensor factorizations have been previously applied in web link analysis
\cite{KoBaKe05} and also in social networks for the analysis of
chatroom \cite{AcCaYe06} and email communications
\cite{BaBeBr07a,Sun2009}.  In these applications tensor factorizations are
used as exploratory analysis tools and do not address the link
prediction problem.

%%% Local Variables:
%%% mode: latex
%%% TeX-master: "paper"
%%% TeX-command-default: "PDFLaTeX"
%%% End: 
\section{Conclusions}
\label{sec:conclusions}

In this paper, we explore several matrix- and tensor-based
approaches to solving the link prediction problem.
We consider author-conference relationships in bibliometric
data and a simulation indicative of user-service relationships in an online context as example
applications, but the methods presented here also have applications in
other domains such as predicting Internet traffic, flight
reservations, and more.
For the matrix methods, our results indicate that using a temporal model
to combine multiple time slices into a single training matrix is superior to
simple summation of all temporal data. We also show how to extend Katz to bipartite
graphs (e.g., for analyzing relationships between two different types of nodes)
and to efficiently compute an approximation to Katz based on
the truncated SVD.
However, none of the matrix-based methods fully leverages
and exposes the temporal signatures in the data. We present
an alternative: the CP tensor factorizations.
Temporal information can be
incorporated into the CP tensor-based link prediction analysis to gain a perspective
not available when computing using matrix-based approaches.

We have considered these methods in terms of their AUC scores and the
number of correct predictions in the top scores. In both
cases, we can see that all the methods do quite well on the DBLP data set. Katz has the
best AUC but is not computationally tractable for large-scale
problems; however, the other methods are not far behind.
Moreover, TKatz-CWT is best for predicting new links in the DBLP data.
Our numerical results also show that the tensor-based
methods are competitive with the matrix-based methods in terms of
link prediction performance.

The advantage of tensor-based methods is that they can better capture
and exploit temporal patterns. This is illustrated in the user-service
example. In this case, we accurately predicted
links several days out even though the underlying dynamics of
the process was much more complicated than in the DBLP case.

The current drawback of the tensor-based
approach is that there is typically a higher computational cost
incurred, but the software for these methods is quite new and will no
doubt be improved in the near future.

%%% Local Variables:
%%% mode: latex
%%% TeX-master: "paper"
%%% TeX-command-default: "PDFLaTeX"
%%% End:    % Conclusions
\section*{Acknowledgments}
  We would like to thank the anonymous referees who provided many insightful comments that helped improve this manuscript.

  This work was funded by the Laboratory Directed Research \&
  Development (LDRD) program at Sandia National Laboratories, a
  multiprogram laboratory operated by Sandia Corporation, a Lockheed
  Martin Company, for the United States Department of Energy's
  National Nuclear Security Administration under Contract
  DE-AC04-94AL85000.  

%%% Local Variables:
%%% mode: latex
%%% TeX-master: "paper"
%%% TeX-command-default: "PDFLaTeX"
%%% LaTeX-command-style: (("." "pdflatex"))
%%% End: 

%% ----------------------------------------------------------------------
%% Bibliography
%% ----------------------------------------------------------------------

\end{document}